\documentclass[reqno,11pt]{amsart}
\usepackage{amsmath, latexsym, amsfonts, amssymb, amsthm, amscd}

\usepackage{xcolor}
\usepackage{enumitem}
\usepackage{graphics,epsf,psfrag}
\usepackage{graphicx}
\usepackage{hyperref}
\usepackage{todonotes}

\numberwithin{equation}{section}

\newtheorem{theorem}{Theorem}[section]
\newtheorem{lemma}[theorem]{Lemma}
\newtheorem{proposition}[theorem]{Proposition}
\newtheorem{cor}[theorem]{Corollary}
\newtheorem{rem}[theorem]{Remark}

\newcommand{\ind}{\mathbf{1}}

\renewcommand{\tilde}{\widetilde}

\newcommand{\cC}{{\ensuremath{\mathcal C}} }

\newcommand{\cL}{{\ensuremath{\mathcal L}} }

\DeclareMathSymbol{\leqslant}{\mathalpha}{AMSa}{"36} 
\DeclareMathSymbol{\geqslant}{\mathalpha}{AMSa}{"3E} 
\DeclareMathSymbol{\eset}{\mathalpha}{AMSb}{"3F}     
\newcommand{\dd}{\,\text{\rm d}}             

\newcommand{\bbC}{{\ensuremath{\mathbb C}} }

\newcommand{\bbE}{{\ensuremath{\mathbb E}} }

\newcommand{\bbL}{{\ensuremath{\mathbb L}} }

\newcommand{\bbN}{{\ensuremath{\mathbb N}} }

\newcommand{\bbP}{{\ensuremath{\mathbb P}} }

\newcommand{\bbR}{{\ensuremath{\mathbb R}} }

\newcommand{\ga}{\alpha}
\newcommand{\gb}{\beta}
\newcommand{\gd}{\delta}
\newcommand{\gep}{\varepsilon}       
\newcommand{\gp}{\varphi}

\newcommand{\gz}{\zeta}
\newcommand{\gG}{\Gamma}

\newcommand{\go}{\omega}

\newcommand{\gs}{\sigma}

\makeatletter
\def\captionfont@{\footnotesize}
\def\captionheadfont@{\scshape}

\long\def\@makecaption#1#2{%
  \vspace{2mm}
  \setbox\@tempboxa\vbox{\color@setgroup
    \advance\hsize-6pc\noindent
    \captionfont@\captionheadfont@#1\@xp\@ifnotempty\@xp
        {\@cdr#2\@nil}{.\captionfont@\upshape\enspace#2}%
    \unskip\kern-6pc\par
    \global\setbox\@ne\lastbox\color@endgroup}%
  \ifhbox\@ne 
    \setbox\@ne\hbox{\unhbox\@ne\unskip\unskip\unpenalty\unkern}%
  \fi
  \ifdim\wd\@tempboxa=\z@ 
    \setbox\@ne\hbox to\columnwidth{\hss\kern-6pc\box\@ne\hss}%
  \else 
    \setbox\@ne\vbox{\unvbox\@tempboxa\parskip\z@skip
        \noindent\unhbox\@ne\advance\hsize-6pc\par}%
\fi
  \ifnum\@tempcnta<64 
    \addvspace\abovecaptionskip
    \moveright 3pc\box\@ne
  \else 
    \moveright 3pc\box\@ne
    \nobreak
    \vskip\belowcaptionskip
  \fi
\relax
}
\makeatother
\def\writefig#1 #2 #3 {\rlap{\kern #1 truecm
\raise #2 truecm \hbox{#3}}}

\newcommand{\logZ}{\mathtt{z}}

\newcommand{\lar}{\vartriangleleft}
\newcommand{\rar}{\vartriangleright}

\newcommand{\newnorm}[1]{{\left\vert\kern-0.25ex\left\vert\kern-0.25ex\left\vert #1 
    \right\vert\kern-0.25ex\right\vert\kern-0.25ex\right\vert}}

\begin{document}

\title[Ising transfer matrix with balanced disorder]{Lyapunov exponent for products of random Ising transfer matrices: the balanced disorder case}

\author[G. Giacomin and R. L. Greenblatt]{Giambattista Giacomin and Rafael L. Greenblatt}
\address[GG]{Universit\'e de Paris,  Laboratoire de Probabilit{\'e}s, Statistiques  et Mod\'elisation, UMR 8001,
            F-75205 Paris, France}
    \address[RLG]{Scuola Internazionale Superiore di Studi Avanzati (SISSA), Mathematics Area, 
	via Bonomea 265,
 34136 Trieste, Italy}

\begin{abstract}
We analyze  the top Lyapunov exponent of the product of  sequences of two by two matrices that appears in the analysis of several statistical mechanics models with disorder: for example these matrices are    the transfer matrices for
 the  nearest neighbor Ising chain with random external field, and  the free energy density of this Ising chain is the  Lyapunov exponent we consider. 
 We obtain the  sharp behavior of this   exponent in  the large interaction limit when 
the external field is centered:   
this \emph{balanced} case turns out to be   \emph{critical} in many respects. 
 From a mathematical standpoint we precisely identify the  behavior of the top Lyapunov exponent
of a product of two dimensional random matrices close to a diagonal random matrix for which  top and bottom Lyapunov exponents coincide. In particular, the
Lyapunov exponent is only $\log$-H\"older continuous. 
 
\bigskip

\noindent  \emph{AMS  subject classification (2010 MSC)}:
60K37,  
82B44, 
60K35, 
82B27  

\end{abstract}

\maketitle

\section{Introduction and results}

\subsection{Background}
The simple two dimensional matrix of the form  
\begin{equation}
\label{eq:keyM}
M\, =\,M(\gep, Z)\, :=\, 
\begin{pmatrix}
1& \gep\\ 
\gep Z  & Z
\end{pmatrix}\, , 
\end{equation}
with $\vert \gep \vert < 1$ and $Z\ge 0$ repeatedly appears in physics, notably in the Ising model context, with and without disorder.
In particular, it is the transfer matrix of the Ising chain 
(see for example \cite{cf:FV,cf:MWbook}, but we refer to Section~\ref{sec:Ising} for an extended discussion, also about the two dimensional Ising model and the one dimensional quantum case).
In the context of disordered models, $Z$ is a non negative random variable:  we consider the case $(Z_j)_{j=1,2, \ldots}$ IID, so $(M(\gep, Z_j))_{j=1,2, \ldots}$ is a sequence of IID random matrices.
From now on we denote $Z$ a random variable with the same law of the $Z_j$'s and our results rely on the hypotheses that $Z$ has a density and that $\bbE[Z^x]< \infty$ for $x$ in a neighborhood of the origin. It is useful to stick to this framework from now, even if occasionally in this introduction we will take some freedom.

Our work focuses on the top Lyapunov exponent of the product of the random matrices $(M(\gep, Z_j))_{j=1,2, \ldots}$:
\begin{equation}
\label{eq:L_Z}
	\cL_Z(\gep)\, :=\, \lim_{n \to \infty} \frac 1 n \bbE \log \left \Vert M(\gep , Z_1) M(\gep , Z_2) \cdots  M(\gep , Z_n)\right \Vert\, ,
\end{equation}
where $\Vert \cdot \Vert$ is an arbitrary matrix norm. The existence of the limit in \eqref{eq:L_Z} and the independence of the choice of the norm holds under very mild assumptions \cite{cf:BL}: in particular this holds under our hypotheses (of course the case $\gep=0$ can be analyzed by elementary methods). Moreover 
we can pass from $\gep$ to $-\gep$ by conjugation via the diagonal matrix with $(1,-1)$ on the diagonal, so 
$\cL_Z(\gep)=\cL_Z(-\gep)$ and it suffices to  analyze the case  $\gep > 0$.

 In \cite{cf:DH}  B.~Derrida and H.~Hilhorst -- henceforth DH --  
one finds the claim that when
$\bbE[ \log Z]< 0$,  $\bbE [Z]>1$ and $\bbE[Z]<\infty$
\begin{equation}
\label{eq:DH}
\cL_Z(\gep)  \stackrel{\gep \to 0}\sim C_Z \vert \gep\vert ^{2 \ga}\, ,
\end{equation}
where $\ga \in (0,1)$  is the only nonzero real solution of $\bbE[Z^\ga]=1$  and $C_Z$ is a positive constant
(existence and uniqueness is an elementary convexity issue, see however
Remark~\ref{rem:alpha} below). Note that \eqref{eq:DH} means that $\cL_Z(\gep)$ is singular at the origin, while it is real analytic for $\gep \in (0,1)$ (analyticity  follows by applying the main result in \cite{cf:Ruelle}, see also  \cite{cf:CGG} and references therein). For the physical relevance of this singularity we refer to \cite{cf:DH} and to Section~\ref{sec:Ising} below. Mathematically 
instead
the lack of regularity may be seen as a result  of the qualitatively different dynamical properties of the action of  the diagonal $\gep=0$ matrices, all of which have two nontrivial invariant subspaces in common (the two coordinate axes), compared to the the $\varepsilon > 0$ case for which the matrices have no nontrivial invariant subspaces in common (for nontrivial distributions of $Z$).

\smallskip

\begin{rem}
\label{rem:alpha<0}
Since we can pass from $M(\varepsilon,Z)$  to $Z M(\varepsilon,1/Z)$ by conjugation with the matrix with $(1,1)$ on the anti-diagonal and $(0,0)$ on the diagonal, in the case   $\bbE[\log Z] > 0 $, $\bbE [1/Z] > 1$ and $\bbE [1/Z] < \infty$ \eqref{eq:DH} yields 
\begin{equation}
\label{eq:DH<}
\cL_Z(\gep) - \bbE[ \log Z] \stackrel{\gep \to 0}\sim C_{1/Z} \vert\gep\vert ^{-2 \ga}\, ,
\end{equation}
where $\ga \in (-1,0)$  is the only nonzero real solution of $\bbE[Z^\ga]=1$.
\end{rem}

\smallskip

The claim in \cite{cf:DH} goes beyond the case of $Z$ variables with a density. Possibly \eqref{eq:DH} holds just assuming that 
$Z$ is not restricted to a discrete set of values of the form $\left\{ 0 \right\} \cup \left\{ \dots,z^{-1},1,z,z^2,\dots \right\}$. In this direction, 
in \cite{cf:DH}  it is pointed out  that the  case in which $Z$ takes just two values $0$ and $z>0$ is explicitly solvable. In this case  $C_Z$ must be replaced by  a log-periodic function of $\gep$. 

  In \cite{cf:GGG}, \eqref{eq:DH} has been proven under the assumption that $Z$ has a compact support bounded away from zero and that $Z$ has a $C^1$ density. 
  Beyond  the several motivations set forth in \cite{cf:CPV,cf:DH} (and partly presented in Section~\ref{sec:Ising} below), this result identifies a singular behavior of the Lyapunov exponents and enters the field of inquiry into the regularity of Lyapunov exponents (e.g.\ \cite{cf:lepage,cf:BaDu}) as a case in which the singularity is sharply identified. 
  
\smallskip

\begin{rem}
\label{rem:alpha}
The function $x \longrightarrow \bbE[Z^x]$ is convex and, by  assumption, it is bounded in a neighborhood of the origin.  It takes value one at the origin where its derivative  is equal to $\bbE[\log Z]$. Therefore if
$\bbE[Z]< \infty$ then  the function is continuous in the whole interval $(0,1)$, so 
$\bbE[ \log Z]< 0$ and $\bbE [Z]>1$ do imply that $\bbE[Z^x]=1$ has only one solution in $(0,1)$. 
 The issue of whether \eqref{eq:DH} holds also when $\bbE[Z]=\infty$ but $\bbE[Z^x]=1$ has a solution in $(0,1)$ (which is unique by convexity)
 may be approachable with the techniques used in \cite{cf:GGG}, which however is restricted to compactly supported random variables and removing this assumption up to allowing $\bbE [Z]=\infty$ does not appear to be straightforward.
 \end{rem}

\smallskip

The results presented up to now are for choices of $Z$ such that $\bbE[\log Z]\neq 0$ and such that $\vert \ga \vert <1$ (we will discuss in Section~\ref{sec:overview} what is known about the case $\vert \ga \vert\ge 1$).  In this work we focus on the case $\bbE[\log Z]=0$ and we will refer to this case, in agreement with the physical literature (see e.g. \cite[p.~1220]{cf:NL}), as the $\ga=0$ case:  in this case   
$\bbE[Z^x]=1$ implies $x=0$ and  the non trivial solution matches  the trivial one. This is a first reason to consider the case as \emph{critical}. We will see that the random matrix product we consider is 
  is directly linked to random walk in random environment or random iteration problems  \cite{cf:BBE,cf:BBD,cf:BDMbook,cf:Sinai}, see also Sections~\ref{sec:firstlook} and \ref{sec:strategy}:   the case $\bbE[\log Z]=0$ has a critical nature in this context and and it is sometimes   referred to as \emph{balanced}.   The critical nature of the balanced case in the Ising context is discussed in Section~\ref{sec:Ising}.

In view of \eqref{eq:DH}, one definitely expects  that $\cL_Z(\gep)$ is singular when  $\gep \to 0$ also if $\bbE[\log Z] =0$. 
In  this work we obtain a mathematical control on the  behavior of  $\cL_Z(\gep)$ for $\gep \searrow 0$ for a wide class of $Z$ such that $\bbE[\log Z]=0$: in particular, we capture the sharp leading asymptotic behavior of the singularity. 

\subsection{The main result}
\label{sec:maires}


For what follows it is  more natural to work with $\logZ_j:= \log Z_j$ and let $\logZ$ be a variable that is distributed like the $\logZ_j$'s.   On $\logZ$ we assume: 
\smallskip

\begin{enumerate}[leftmargin=*]
\item  exponential integrability, namely that there exists  $\gd>0$  such that 
\begin{equation}
\label{eq:LR}
\bbP\left( \left \vert  \logZ\right \vert > x\right) \,=\, O \left( \exp\left( -\gd x\right)\right) \, .
\end{equation} 
\item $\logZ$ has a density  $\gz$ which  is uniformly $\theta$-H\"older continuous, for a $\theta \in (0, 1]$, i.e.\ 
\begin{equation} 
\label{eq:Holder}
\sup _{x\not= y}
\frac{\vert \gz(x)-\gz(y)\vert}{\vert x - y \vert ^\theta} \,<\,  \infty \, .
\end{equation}  
\end{enumerate}
\smallskip

We remark that these hypotheses imply  that $\gz$ is locally bounded and that $\gz(x)=O(\exp(-\gd' \theta\vert x \vert/(1+ \theta))$ for every $\gd'\in (0, \gd)$ (see Lemma~\ref{th:Leonardo}(1)). In particular,  $\gz \in \bbL^p$ for every $p\in [1, \infty]$.  

\smallskip

Of course we assume $\bbE[\log Z]=0$, that is  $\bbE[\logZ]=0$.

\medskip

\begin{theorem}
\label{th:main}
There exist three constants $\kappa_1>0$, $\kappa_2 \in \bbR$  
and $\eta \in (0,1)$ such that, for $\gep \to 0$,
\begin{equation}
\label{eq:main}
\cL_Z(\gep)\, =\,  \frac {\kappa_1}{  \log (1/ \vert\gep\vert ) +\kappa_2 } + O\left( \vert\gep\vert^\eta\right)\, .
\end{equation}
\end{theorem} 

\medskip

 The constants $\kappa_1$, $\kappa_2$ and $\eta$ depend only on the law of $Z$.
For $\kappa_j$  we have a semi-explicit expression, see \eqref{eq:constants},
and for $\eta$ we give an explicit lower bound. 

\smallskip

Theorem~\ref{th:main} shows in particular that $\cL_Z(\cdot)$ is not H\"older continuous: this is not in  contrast with
  \cite{cf:BaDu,cf:lepage} 
  because we are looking at the neighborhood of diagonal random matrices (hence lacking irreducibility: the two axes are invariant) and for in which there is no separation between the two Lyapunov exponents (in the Osedelec sense  \cite[Ch.~IV]{cf:BL}). 
  Theorem~\ref{th:main} actually provides an ensemble of examples for which there is a  sharp control on a $\log$-H\"older singular behavior of Lyapunov exponents:
 with respect to this we mention the bounds in this spirit obtained in the context  of Anderson localization (\cite{cf:CS} and references therein).

In  \cite[(4.34) and pp. 1218-1220]{cf:NL}  one finds the statement \eqref{eq:main}
 for a specific class of laws of $Z$ (superposition of one (bi-)exponential and one Dirac delta) and $\bbE[ \log Z]=0$.
 The expressions for $\kappa_1$ and $\kappa_2$ are more explicit, however we are unable to say whether 
 the computations in \cite{cf:NL} can be made into a mathematical proof (in particular, the laws of $Z$ considered in  \cite{cf:NL}
 do not satisfy our hypotheses). Moreover in \cite{cf:CTT} a class of random matrices for which the Lyapunov exponent can be  expressed explicitly is identified: this includes a model which is close the ones we consider, i.e.\ a model in which $\logZ$ is a bi-exponential, but $\gep$ is multiplied by an independent exponential random variable. The small $\gep$ asymptotic behavior in this case matches \eqref{eq:main}.

 The result  \eqref{eq:main} is  known in the simpler framework of the so called  \emph{weak disorder limit} (see Section~\ref{sec:overview}). 
 Moreover, the \emph{2-scale approach} of \cite{cf:DH} can be generalized and leads to a prediction  in the spirit  of \eqref{eq:main} for \emph{arbitrary} distributions of $Z$ such that $\bbE[\log Z]=0$  \cite{cf:Derrida-private}: this is very relevant to us because, as we explain next and above all in Section~\ref{sec:strategy}, 
 our proof of  \eqref{eq:main} stems from 
 the DH 2-scale idea.

 \subsection{A first look at the strategy of proof}
 \label{sec:firstlook}

 As we just announced, like  in   \cite{cf:GGG} we  start from  the 2-scale idea of \cite{cf:DH}. What is done in \cite{cf:DH} for $\bbE[\log Z] \neq 0$ is to guess a probability measure -- we call it the DH probability -- that should be sufficiently close to the Furstenberg probability (on the projective space, i.e.\ simply the sector  $[0, \pi/2]$ because we work with positive matrices in dimension two) with which one can explicitly express the Lyapunov exponent \cite{cf:BL}. The measure essentially concentrates near $0$ (we are assuming $\bbE[\log Z] < 0$, it would be near $\pi/2$ if $\bbE[\log Z] > 0$), but a much finer description is needed. The 2-scale idea is about gluing together   two $\gep \searrow 0$ limit problems: one  near $0$ and another  near $\pi/2$. Both problems are non trivial and they require a quantitative understanding of the invariant measure of 
a chain that appears in  the context of one dimensional  Random Walks in Random Environment (RWRE) \cite{cf:decalanetal,cf:kozlov,cf:KKS}  and in  random affine iterations \cite{cf:BDMbook}. Even if the problems at $0$ and $\pi/2$ are in a sense \emph{dual problems}, when $\bbE[\log Z]\neq 0$ they are  different in nature: 
one of the two chains is positive recurrent, the other is transient.
  In \cite{cf:GGG} we first gave a rigorous construction of the  DH probability -- this is essentially an asymptotic matching problem --
    and then we showed that this probability is \emph{sufficiently close} to the Furstenberg probability to yield 
    \eqref{eq:DH}. 
In  the key second step we exploited a contraction property, under the random matrix action, of a suitable norm that depends on a parameter $\gb\in (0,\ga)$:
the contraction factor is precisely $\bbE[Z^\gb]<1$. 
\medskip

\begin{figure}[h]
\centering
\includegraphics[width=13.5 cm]{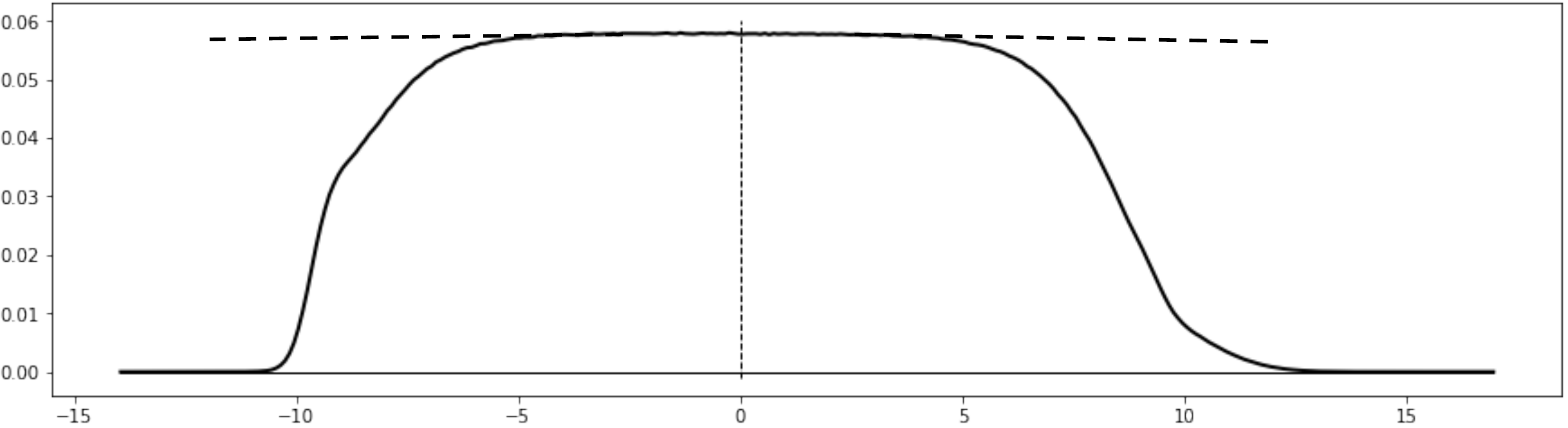}
\caption{\label{fig:2} 
The solid line gives density of the invariant probability for 
$k=10$ and $\logZ= \xi (-1/2+3\eta/10)+(1-\xi)(\eta+1)$, with $\xi$ a Bernoulli with success probability $2/3$ and   $\eta$ a standard Gaussian, with $\xi$ and $\eta$ independent. The distribution of $\logZ$ is asymmetric and bimodal. 
The densities of the two non normalizable measures that build the DH measure are (essentially) given by the prolongation of the plateau to the right for one measure, and to the left for the other one (dashed lines). 
}
\end{figure}

When $\bbE[\log Z]=0$ there are two major changes with respect to the $\bbE[\log Z]\neq 0$ case :
\smallskip

\begin{enumerate} [leftmargin=* ]
	\item The two limit problems change nature: they become two qualitatively identical problems -- if $Z$ has the same law as $1/Z$ they are the same problem -- and they are  both null recurrent. These chains are associated to a balanced RWRE, also known as \emph{Sinai walk}, and they are a particular critical random difference equation \cite{cf:BBE,cf:BBD}: 
  these chains do  not have an invariant probability, but each of them  does have a \emph{unique} $\gs$-finite invariant measure on which we need a sharp control in order to perform the gluing procedure.  The gluing procedure builds the 
  DH measure, which is a probability measure and is close to the invariant probability of the chain we are interested in (which is positive recurrent), from the two $\gs$-finite invariant measures, see Fig.~\ref{fig:2} .  
\item $\bbE[\log Z]=0$ means  $\ga=0$, so there is no room to apply the contractive procedure of \cite{cf:GGG}.  As a matter of  fact,  this 
lack of contraction makes clear one of the reasons to consider the 
problem as critical (and this is why it is considered critical also in \cite{cf:BBE,cf:BBD}). 
\end{enumerate}

\smallskip
We will nevertheless take the DH path: we are going to explain in  \S~\ref{sec:strategy} how we do it, notably how we deal with the lack of the contraction property.

\subsection{More on the approach: the DH strategy for the balanced case}
\label{sec:strategy}
We now go into some of the details of the construction of the DH probability and we explain the idea of the proof that this measure is \emph{sufficiently close} to the Furstenberg probability. 

 A change of perspective on the problem is notationally useful: we set 
\begin{equation}
k\, :=\, -\log \gep\, ,
\end{equation}
so that 
  the original matrix $M$ of \eqref{eq:keyM}
becomes
 \begin{equation}
\label{eq:matrix2}
\begin{pmatrix}
1& \exp(-k)\\ 
\exp(-k+\logZ)  & \exp(\logZ)
\end{pmatrix}\, ,
\end{equation}
and if we parametrize $P(\bbR^2)\cong \bbR$ with the coordinates $(1, \exp(x))^{\mathtt{T}}$ we readily see that
the action of the matrix \eqref{eq:matrix2} is 
\begin{equation}
\label{eq:map}
x \, \mapsto \, z+ \log \left(
\frac{e^{-k}+e^x}{1+ e^{x-k}}
\right)\, =\, z+ h_k(x)\, .
\end{equation}
We observe that $h_k(\cdot)$ is odd and 
$x \mapsto x- h_k(x)$, which  is also odd,  is small if $x \in (-k,k)$ and far from the boundary points $\pm k$:
\begin{equation}
\label{eq:almost-id}
x-h_k(x)\, =\,  \log \left(\frac{1+ e^{x-k}}
{1+e^{-x-k}}
\right)\, 
=\,  O\left( e^{-(k-x)}\right)+ O\left( e^{-(k+x)}\right)\, ,
\end{equation}
so $h_k(x)$ is \emph{very close} to $x$ on an interval that approaches $\bbR$ in the limit $k \to \infty$ (see Fig.~\ref{fig:1}).

\medskip

Denote by $X=(X_n)_{n=0,1, \ldots}$ the Markov chain generated by 
the map \eqref{eq:map}, that is 
\begin{equation}
\label{eq:XMC}
X_{n+1}\, =\,  \logZ_{n+1} + h_k\left(X_n \right)\, .
\end{equation}
Since the image of $h_k$ is $(-k,k)$ the process $(X_n)$ \emph{hardly leaves} $(-k,k)$ and when the process is in $(-k,k)$ and far from the boundaries it is \emph{close to being} a random walk with centered increments. 


$X$ is  an irreducible positive recurrent Markov chain (see the beginning of Appendix~\ref{sec:various}) and via
its invariant probability $\nu_k$ one has that
\begin{equation}
\label{eq:Lyap-k}
 \cL(k)
 :=
\cL_Z(\exp(-k)))\,=\, \int _\bbR \log(1+ \exp(-k - x)) \nu_k(\dd x) \, .
\end{equation} 
 The formula \eqref{eq:Lyap-k} will be further explained later on (see in particular \eqref{eq:L_k}), but it is what follows by specializing the Furstenberg formula for the top Lyapunov exponent \cite[Th.~3.6]{cf:BL} to our context.

Since $X$ is  close to being a symmetric random walk in the bulk and because of the strong containing effect at the boundary, it is natural to expect that the invariant probability is going to be close to the Lebesgue measure times a suitable constant $1/C_k$ in the bulk, and this measure should decay quickly outside  $(-k, k)$. If this is the case $C_k \sim 
2k$ in the $k\to \infty$ limit.  If we insert this guess into \eqref{eq:Lyap-k} 
we readily see that the leading contribution comes from $x$ close to $-k$ on the scale $k$: by this we mean an interval centered in $-k$ of diverging length $o(k)$.  
We are therefore interested in focusing on how the process looks from $-k$. So we consider 
$(X_n+k)$ and  we readily see that  the process that appears for  $k\to \infty$ is 
\begin{equation}
\label{eq:iterY-0}
Y_{n+1}\, =\,  \logZ_{n+1} + h\left(Y_n \right)\ \ \ \text{ with } h(y)\, =\, y+ \log(1+\exp(-y))\, .
\end{equation}
This new Markov chain, of which we will give a detailed treatment, has a strong repulsion at zero, forcing $Y$ to live most of time in the positive semi-axis. But there is no mechanism that bounds $Y$ from above: in fact, $Y$ is null recurrent and the non normalizable invariant measure does approach a multiple of the Lebesgue measure far from the origin.  The
random iteration \eqref{eq:iterY-0} is the critical
  random difference  equation that emerges in the analysis of 
Sinai walks \cite{cf:BBE,cf:BBD,cf:BDMbook,cf:Sinai}; consequently, the literature is extensive, notably  on the behavior of the invariant measure of the $Y$ process at $+\infty$. However the focus for us is twofold: 
\smallskip

\begin{enumerate}[leftmargin=*]
\item characterizing the local part of the invariant measure and the behavior at $-\infty$;
\item obtaining a sharp estimate on the behavior at $+\infty$. 
\end{enumerate} 

\smallskip

While point (1) is central because it determines the leading behavior of \eqref{eq:Lyap-k}, point (2) is as important because 
the invariant measure of $Y$ just provides  the DH guess on the negative and positive  semi-axes, separately (they are essentially symmetric problems). 
But we need a 
 sharp asymptotic analysis at $+\infty$ of the two measures, i.e.\ at the origin which is very far from both $+k$ and $-k$, to glue them together. The results available on this problem are obtained in too general a context and they are too weak for our purposes: we therefore perform  an \emph{ad hoc} analysis. 
 
Once the  DH probability  $\nu^{\textrm{DH}}_k$ (explicit!) is built, its closeness  to the invariant probability $\nu_k$ (not explicit!) has to be established. What is directly accessible is the action of the (one step) Markov Kernel $T=T_k$ on $\nu^{\textrm{DH}}_k$ and we certainly want $T \nu^{\textrm{DH}}_k-\nu^{\textrm{DH}}_k$ \emph{small}. In \cite{cf:GGG}  this closeness is estimated for $\ga>0$ in terms of a 
norm $\newnorm{\cdot}_\gb$ indexed by $\gb\in (0,  \ga )$. And the key point is that $\newnorm{\nu^{\textrm{DH}}_k-\nu_k}_\gb$ is bounded above by $c_\gb \newnorm{T \nu^{\textrm{DH}}_k-\nu^{\textrm{DH}}_k}_\gb$, with
 $c_\gb=1/(1-\bbE[Z^\gb])_+$ (if $\ga <0$, i.e.\ $\bbE[Z^\ga ]>1$,  we  work with $\bbE[Z^{-\gb}]$). 
 $\newnorm{\cdot}_0$ is well defined too: 
it is simply  the $\bbL^1$ norm of the primitive of $\cdot$ which is in $\bbL^1$. And in  fact $\newnorm{\cdot}_0$
 is, or could be, a good norm for our purposes: the problem is that   $c_0=+\infty$.

To get around this problem we take an approach that avoids using contraction properties. In fact we show that  the bound holds with $c_0=c_0(k)$ which is $O(k^2)$ (up to logarithmic corrections). The divergence of $c_0(k)$ can be overcome if $\newnorm{T \nu^{\textrm{DH}}_k-\nu^{\textrm{DH}}_k}_0$ decays faster that $k^{-2}$. With our hypotheses, this decay is exponential in $k$. 

\smallskip

\begin{rem}
\label{rem:sym}
At this stage the dominant role of $-k$  with respect to $+k$ may appear strange. But this is just an artifact of the choice of 
\eqref{eq:Lyap-k} which is the formula that one obtains when looking at the exponential growth of the $(1,1)$ entry of the matrix.  But $+k$ takes the leading stand if we consider the formula stemming out of  the exponential growth of the $(2,2)$ entry, see  \eqref{eq:L_k}. Of course, the Lyapunov exponent does not depend on the choice of entry. This is discussed further in  Remark~\ref{rem:lyapsim}.
\end{rem}

\subsection{Exactly solvable Ising models and matrix products}
\label{sec:Ising}

In \cite{cf:CPV,cf:DH} the role of the two by two matrix \eqref{eq:keyM} in the solution or analysis of several physical models is discussed and exploited. 
Here we focus on cases related to  the Ising model. 

\begin{enumerate}[start=1,label={(I\arabic*)},leftmargin=*]
\item \label{Is1} If $d=1$  the transfer matrix of the  Ising model with external magnetic field $h$ and nearest neighbor interaction potential $J$ is equal to  $e^{h+J}$ times \eqref{eq:keyM}, with $\gep= e^{-2J}$ and $Z=e^{-2h}$. Therefore, 
without disorder the free energy density of the model in the thermodynamic limit is equal to $h+J= (-\log \gep + \log Z)/2$ plus
\begin{equation}
	\lim_{n \to \infty}
	\frac{1}{n}
	\log
	\mathrm{Trace}
	\left[ \left(M\left( \varepsilon, Z \right) \right) ^n\right]
	\, ,
\end{equation}
which in turn is given by the largest among the two eigenvalues of $M( \varepsilon, Z)$.

With disorder in the magnetic field, $Z$ takes a sequence of IID values $Z_1,Z_2,\dots $ and the free energy density  coincides with $(-\log \gep + \bbE[\log Z])/2+\cL_Z(\gep)$. Therefore the $\gep\searrow 0$ analysis captures the behavior of the free energy density of the one dimensional Ising chain with disordered magnetic field in the infinite coupling interaction limit.
\item	\label{Is2} For the Ising model in two dimensions (on the square lattice) with horizontal and vertical couplings $J_1$ and $J_2$ the transfer matrix is much larger, but for $h=0$ it can be rewritten in terms of a tensor product of matrices of the form \eqref{eq:keyM} with different values of $\gep$: explicit eigenvalue computations lead to the celebrated Onsager determination of the critical behavior of the model. 
Remarkably, this tensor product structure stands also when  the couplings vary in one of the coordinate directions but remain constant in the other (so-called \emph{columnar disorder}; this model was first understood in \cite{cf:MW} and is now known as the McCoy-Wu model \cite{cf:MW,cf:MWbook,cf:SM}). But of course, the problem is no longer solvable by computing eigenvalues: now the free energy density can be expressed in terms of Lyapunov exponents. More precisely, 
	in the non disordered case the free energy density in the thermodynamic limit is given by an integral over a parameter $\theta$ of the largest eigenvalue of $M( \varepsilon, Z)$, with an adequate choice of $\gep(\theta, J_1)$ and $Z=Z(J_1, J_2)$.
	$\theta$ can be understood as a Fourier transform parameter, and long distances (and hence critical phenomena) appear in the behavior of this eigenvalue for $\theta \to 0$, and  $\lim_{\theta \to 0}\gep(\theta, J_1)=0$. 
	In the disordered case, the formula is formally the same, with largest eigenvalue replaced by top Lyapunov exponent 
	 (several details can be found in \cite[Appendix~A]{cf:CGG}, but of course also in \cite{cf:MW,cf:MWbook,cf:SM}). 
\item 	\label{Is3}
	The partition function of the  one-dimensional  transverse field  quantum Ising model can be written using the Trotter product formula as a limit of that of the two-dimensional classical Ising model. Disorder in the one-dimensional quantum model, either in the transverse field or in the couplings, can be shown to correspond to columnar disorder in the two-dimensional classical model.
	So, once again the same matrix product appears and the
	 the analysis of the ground state, i.e. a suitable zero temperature limit,  of the  quantum model corresponds to analyzing the $\gep\to 0$ limit. We refer to 
	 \cite{cf:Luck,cf:Fisher} for several details on the correspondence between quantum chain and two dimensional model. 
\end{enumerate}

\smallskip 

In case  \ref{Is1}), with or without disorder,  there is no phase transition, i.e. the free energy is a real analytic function. This is elementary to establish without disorder. In the disordered case one has to resort to more advanced techniques 
  \cite{cf:Ruelle} and a crucial point is the (strict!) positivity of all the entries of the two by two random matrices. In fact, positivity is lost precisely 
   when the coupling strength diverges and the emerging singularity (that we study) may be seen as a pseudo-critical behavior. 


 On the other hand  it is well known that, without disorder,  for $d=2$ and  in the $d=1$ quantum case (see \ref{Is2} and \ref{Is3})   there is a phase transition. This is possible due to the fact that $\gep(\theta, J_1)=0$ for $\theta=0$, 
 so we are no longer in the context of matrices with positive entries and one verifies by an explicit computation that the transition happens if and only if 
 at $\gep=0$ the two eigenvalues coincide. The important claim by McCoy and Wu  \cite{cf:MW,cf:MWbook} (see also the \emph{weak disorder limit} in Section~\ref{sec:overview}) is that 
 also in the columnar disorder the transition happens if and only if at $\gep=0$ the two Lyapunov exponents coincide, and this means precisely $\bbE[\log Z]=0$. Mathematically this claim by McCoy and Wu is an important open problem. We do analyze exactly this critical case, but we do not study  the (non-)analytic dependence of the Lyapunov exponent on the inverse temperature $\gb$. In fact there is no $\gb$ in our model because it is mathematically irrelevant: $\gb$ appears in the model as a multiplicative factor on the interactions (i.e., 
 $Z=Z(\gb J_1, \gb J_2)$). Since we are looking at the singular behavior in $\gep$, there is no loss of generality in setting $\gb=1$. 

\subsection{Generalizations, more on the literature and organization of the paper}
\label{sec:overview}

\subsubsection{Distributions with $|\alpha| \ge 1$}
When $\bbE[\log Z] < 0$  but $\bbE[Z] \le 1$ then  $\bbE[Z^\ga]=1$ may have a (unique) solution  solution for $ \ga \ge 1$. Existence 
requires integrability conditions (see Remark~\ref{rem:alpha}) and  also that  $\bbP(Z >1)>0$. The situation is analogous for 
$\bbE[\log Z] > 0$ and $\bbE[1/Z] \le 1$.
One can find in \cite{cf:DH} an argument  in favor of  $\cL_Z(\gep)= \sum_{j=0}^{\lfloor \ga  \rfloor} C_j \gep^{2j}+ C_\ga \gep^{2 \ga } +o( \gep^{2 \ga })$ for $\gep \searrow 0$, at least for $\ga  \not \in \bbN$ (and analogous formula for $\ga <0$).
 Results in this direction can be found in \cite{cf:benjamin}. 
 Let us point out, however, that identifying the singular behavior, that is capturing the term $C_\ga \gep^{2 \vert \ga \vert}$, is an open problem and the approach we employ does not seem to be appropriate when this is a subleading term (i.e., $\vert \ga \vert \ge 1$). 

\subsubsection{Weak disorder limits}
The weak disorder limit  corresponds in dynamical terms to a very slow dynamics. This allows to rescale time and the arising dynamics is a two dimensional system of stochastic ODEs that can be  solved. This limit dynamics has the remarkable property that the Furstenberg probability has an explicit formula which leads to  an expression for the Lyapunov exponent  in terms of Bessel functions (to our knowledge this appeared first in the works of McCoy and Wu \cite{cf:MW}, but similar computations for related systems were already in the literature, 
see\cite{cf:CGG,cf:CLTT,cf:CPV,cf:Luck,cf:Sadel} for recent works on weak disorder limits and historical overviews). Deriving the full asymptotic expansion of the Lyapunov exponent is then a somewhat cumbersome, but straightforward, exercise \cite[Prop.~1.3]{cf:CGG}. Particularly relevant for us is the $\gep\searrow 0$ limit behavior of the Lyapunov exponent for $\ga=0$  \cite[(1.11)]{cf:CGG} which can be refined to \cite[end of Sec.~4]{cf:CGG} 
\begin{equation}
\label{eq:wd}
  \frac {1}{ 4\left(\log (1/ \gep)  - \log 2 - \gamma\right)} + O\left( \gep^2 \right)\,,
\end{equation}
where $\gamma$ is the Euler-Mascheroni constant and the $O(\gep^2)$, which can be explicitly expressed, is non zero.  
Of course the last formula should be compared to our main result \eqref{eq:main}.

As already pointed out in \cite{cf:CGG}, the extremely sharp matching of the \emph{finite disorder} case and the \emph{infinitesimal disorder} case is mathematically quite surprising: not only the leading order matches, but the structure of the subleading terms is the same (to the extent of what is known). However there does not appear to be any way  to recover finite disorder results from infinitesimal disorder computations. This is particularly unfortunate for the $d=2$ columnar disorder case  (see  \cite{cf:MW,cf:MWbook,cf:SM} and Sec.~\ref{sec:Ising}): McCoy and Wu used the weak disorder limit to infer results about the model and setting forth a precise and highly non trivial prediction for the critical behavior that represents a challenge for mathematicians (see \cite[App.~A]{cf:CGG} and references therein). 

\subsubsection{What about lower regularity of the distribution of the disorder?}
Our approach does not capture only the leading behavior, but also all the  subleading corrections in powers of 
$1/\log(1/\gep)$,  much like \eqref{eq:wd} which stems instead from an exact computation. 
Clearly in estimating the error term there seems to be a lot of room: for example, it would be sufficient to show (\emph{the one step estimate}, see Sec.~\ref{sec:strategy})
that $\Vert T \nu^{\textrm{DH}}_k-\nu^{\textrm{DH}}_k\Vert_1 =O(k^{-c})$, for a $c>2$. But our approach naturally yields an exponential estimate, see \eqref{eq:2addenda}. This approach 
(i.e., estimating the inverse  Laplace transform) strongly relies on an appropriate regularity of $\gz$, i.e.\ asymptotic decay in Fourier space (the imaginary direction for the Laplace transform), and  yields leading term, subleading corrections and 
exponential bounds (in the $k$ variable). It is really not clear whether our one step bound procedure could go through with less regularity.

     
\subsubsection{More general matrices}
All we did can be greatly generalized (and this discussion applies to \cite{cf:GGG} too).  Matrices of the form
\begin{equation} 
\begin{pmatrix}
\widetilde Z& \gep \widetilde Z\\ 
\gep Z  & Z
\end{pmatrix} \ \ \ \text{ and } \ \ \
\begin{pmatrix}
1& \gep Z'\\ 
\gep Z''  & Z
\end{pmatrix}\, , 
\end{equation}
can be dealt with in a straightforward way under suitable hypotheses: for the first example, if  
$\widetilde Z>0$ too, with no independence requirement with $Z$,  simply because $\widetilde Z$ can be factored out and the key variable becomes $Z/ \widetilde Z$. For the second example, the analysis does not contain much novelty if 
$Z$, $Z'$ and $Z''$ are independent and positive and if we require on  the marginal laws of $Z'$ and $Z''$ the same integrability and regularity we have required for $Z$. A complete analysis of two dimensional matrices with positive entries approaching for $\gep\searrow 0$ 
a diagonal matrix might be of interest (but it is probably rather cumbersome and the Ising model motivation is lost). Some   higher dimension matrices generalizations  are treated in \cite[App.~A]{cf:benjamin}.

\subsubsection{Organization of what follows}
In Section~\ref{sec:XY} we study the main Markov chain $X$. We state   in Proposition~\ref{th:bound} the crucial bound that we call \emph{reduction to one step estimate}. We also study the auxiliary chain $Y$ that is central in the construction of the DH probability. 

In Section~\ref{sec:DH} we construct the DH probability and perform the one step bound: 
the proof of our main result, Theorem~\ref{th:main}, is in this section (right after Proposition~\ref{th:onestep}). 

In Section~\ref{sec:bound} we prove the \emph{reduction to one step estimate}, i.e.\ Proposition~\ref{th:bound}.

\section{The underlying Markov chains}
\label{sec:XY}

In general, if $Z$ is a random variable, we use $G_Z(x)=\bbP(Z>x)$ and $F_Z(x)=1-G_Z(x)$. This notation is used also for finite  measures $\mu$: 
$G_\mu(x)= \mu((x, \infty))$ and 
$F_\mu(x)= \mu((-\infty, x])$. We will work also with non normalizable ($\gs$-finite) measures, notably with measures $\mu$ such that $\mu((-\infty, x])< \infty$ for every $x$, and the notation $F_\mu(x)$ will be used also in this case.

\subsection{About the main Markov chain $X$}
We start the analysis of the $X$ chain defined by \eqref{eq:XMC}.
\medskip

\begin{lemma}
We have 
\begin{equation}
\label{eq:GXTn}
T G_{X_n}(x)\, :=\,G_{X_{n+1}}(x)\, =\,  G_{\logZ}(k+x) +\int_\bbR G_{X_n}(y)h_k'(y) \gz\left( x-h_k(y)\right) \dd y
\,  ,
\end{equation}
and
\begin{equation}
\label{eq:FXTn}
T F_{X_n}(x)\, :=\,F_{X_{n+1}}(x)\, =\,  F_{\logZ}(-k+x) +\int_\bbR F_{X_n}(y)h_k'(y) \gz\left( x-h_k(y)\right) \dd y
\,  .
\end{equation}
\end{lemma}

\medskip

\noindent
\emph{Proof.} By using that $h_k(\cdot)$ is an increasing bijection from $\bbR$ to $(-k,k)$ we see that
\begin{equation}
\begin{split}
G_{X_{n+1}}(x)\, &=\, \bbP\left( \logZ_{n+1} +h_k\left( X_n\right) >x \right) 
\\
&=\, \bbP \left( x-  \logZ_{n+1}\le -k\right)+\bbP\left(h_k\left( X_n\right) >x-  \logZ_{n+1} ; \, \vert x-  \logZ_{n+1}\vert <k \right) 
\\
&=\,  \bbP \left( \logZ \ge x+k\right)+\bbP\left(X_n >h_k^{-1}\left(x-  \logZ_{n+1}\right) ; \, \vert x-  \logZ_{n+1}\vert <k \right)
\\
&=\,  G_{\logZ}(k+x) + \int_{-k+x}^{k+x} G_{X_n}\left( h_k^{-1}(x-z)\right) \gz(z) \dd z\, ,
\end{split}
\end{equation}
and by the change of variable $z=x-h_k(y)$ we complete the proof of the first identity. The proof of the second identity is of course exactly analogous (or can be directly derived from the first).
\qed

\medskip

The map $T$ can then be extended by linearity to act on $G$ defined by 
$G(x)=G_{\nu_+}(x)- G_{\nu_-}(x)$, $\nu_\pm$ finite measures (so $\vert G(-\infty)\vert < \infty$):
\begin{equation}
\label{eq:GXT}
 T G(x)\, =\, G(-\infty)G_{\logZ}(k+x) + \int_{-k+x}^{k+x} G\left( h_k^{-1}(x-z)\right) \gz(z) \dd z
\, .
\end{equation}
In the applications we consider $\nu$ is  the difference of two probability measures, so $G(-\infty)=0$ and
$T$ reduces to $T_0$:
\begin{equation}
T_0 G (x)\, :=\, \int_\bbR G(y)h_k'(y) \gz\left( x-h_k(y)\right) \dd y\, ,
\end{equation}
which is well defined  also in the slightly different set-up of  $G\in \bbL^1$ because \eqref{eq:Holder} implies that
$\Vert \gz \Vert_\infty < \infty$. 
\smallskip

The following bound will be crucial; Section~\ref{sec:bound} is devoted to its proof.

 \medskip

\begin{proposition}
\label{th:bound}
There exists $C>0$ and $k_0$ such that for every $k \ge k_0$ and 
 $G\in L^1(\bbR; \bbR)$
\begin{equation} 
\label{eq:bound}
\sum_{n=0}^\infty \left \Vert T_0^n G \right\Vert_1 \, \le \, k^2 \left( \log k \right)^C \left \Vert G \right\Vert_1\, .
\end{equation}
\end{proposition} 
 
\medskip

It is well known that the Lyapunov exponent  can be  expressed in terms the invariant probability for the action of the random matrix on the projective circle 
\cite[Ch.~2]{cf:BL} and  for two by two  matrices with positive entries we can work on $(0, \pi/2)$ or, considering the tangent of this angle, 
on $(0, \infty)$ 
(see \cite[(1.6)]{cf:GGG} or \cite[Sec.~2]{cf:DH}). Our parametrization choice  $P(\bbR^2)\cong \bbR$, recall \eqref{eq:matrix2}-\eqref{eq:map}, just corresponds to applying the logarithm to the tangent of the angle and it suffices to apply this change of variables to the expressions in  
\cite{cf:DH,cf:GGG} (which take advantage of the specific form of the matrix under consideration to obtain a simpler expression than the standard Furstenberg formula): we obtain that,
by writing $\nu_k$ for the invariant measure on $\bbR$, the Lyapunov exponent $\cL(k)$ is equal to $L_k[G_{\nu_k}]$ with
\begin{equation}
\label{eq:L_k}
L_k[G]\,:=\, \int_ \bbR \frac{1} {1+e^{k-x}} G(x) \dd x\, =\,  \int_ \bbR \frac{1} {1+e^{k+x}} (1-G(x)) \dd x\,,
\end{equation}
which will be used for 
 $G(x)= \nu((x, \infty))$ and $\nu$ a probability.  
We readily see that 
\begin{equation}
\label{eq:Ly}
\left \vert L_k[G_1] - L_k[G_2 ]\right \vert \, \le \, \Vert G_1 -G_2 \Vert_1\, . 
\end{equation}
We thus have the following important corollary to Proposition~\ref{th:bound}:
\medskip

\begin{cor}
\label{th:bound-cor}
With $C$ and $k_0$ as in Proposition~\ref{th:bound},
for $k \ge k_0$ and any  probability $\gamma $ such that $T G_\gamma - G_\gamma \in \bbL^1$
\begin{equation} 
\label{eq:bound-cor}
\left \vert \cL(k)- L_k[G_\gamma] \right \vert \, \le\, \left \Vert G_{\nu_k} - G _ \gamma \right \Vert _1\, \le  \, k^2 \left( \log k \right)^C \left \Vert T G_\gamma - G_\gamma \right\Vert_1\, .
\end{equation}
\end{cor}

\medskip

\noindent 
\emph{Proof.} 
Since the Markov chain $(X_n)$ is irreducible and positive recurrent (see the beginning of Section~\ref{sec:various}), $\lim _n T^n  G_\gamma (x)=
G_{\nu_k}(x)$ for every $x \in \bbR$ which is a continuity point of $G_{\nu_k}(\cdot)$
 (we will see that 
$G_{\nu_k}(\cdot)$ is continuous, see Remark~\ref{rem:nu}, but at this stage this is not needed since 
the set of discontinuities of $G_{\nu_k}(\cdot)$ is countable, hence of Lebesgue measure zero). 
Therefore, by Fatou's Lemma, we have that 
\begin{multline}
\left \Vert G_{\nu_k} - G _ \gamma \right \Vert _1 \, \le 
\, \liminf _N \left \Vert T^N G_\gamma - G_\gamma \right\Vert_1
\\
 \le \, 
 \sum_{n =0}^\infty  \left \Vert T_0^n \left( T G_\gamma - G_\gamma\right) \right\Vert_1 \, \le \, k^2 \left( \log k \right)^C \left \Vert T G_\gamma - G_\gamma \right\Vert_1\, ,
\end{multline}
and the proof is complete by \eqref{eq:Ly}.
\qed

\medskip


\subsection{Looking from the edge: the reduced chain $Y$}
If we sit on $-k$, that is if we make it our new origin, in the limit as $k\to \infty$ the Markov chain becomes
\begin{equation}
\label{eq:iterY}
Y_{n+1}\, =\,  \logZ_{n+1} + h\left(Y_n \right)\ \ \ \text{ with } h(y)\, =\, y+ \log(1+\exp(-y))\, .
\end{equation}

\begin{figure}[h]
\centering
\includegraphics[width=13.5 cm]{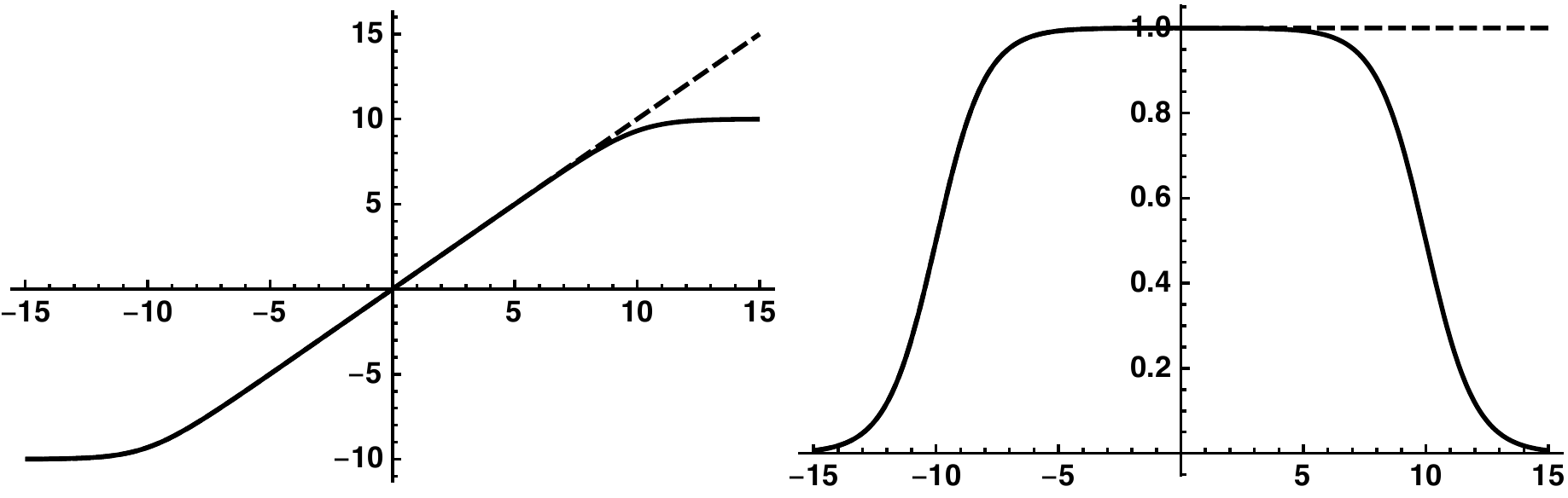}
\caption{\label{fig:1} For $k=10$: on the left the plot of $h_k(\cdot)$ (solid line) and of $x\mapsto h(x+k)-k$ (dashed line); on the right the derivative of the same functions. }
\end{figure}

\medskip

\begin{proposition}
\label{th:Y}
\begin{sloppypar}
	The Markov chain $Y$ has a unique   invariant measure $\nu$ with ${\nu((-\infty,x])< \infty}$ for every $x\in \bbR$, but
	$\nu(\bbR)= \infty$. 
\end{sloppypar}
\end{proposition}
\medskip

For non-normalizable measures uniqueness is of course meant \emph{up to  a multiplicative constant} and we note that $\nu$ is characterized by 
\begin{equation}
\label{eq:fornu}
\int  g(z) \nu (\dd z)\, =\, \int \int g(z +h(y)) \gz(z) \dd z \, \nu (\dd y)\ \ \text{ for every measurable } g \ge 0\,  .
\end{equation}

The proof of Proposition~\ref{th:Y} is  given for completeness in the Appendix: the Markov chain $Y$ is a (very) particular case of the critical   random difference equation problem, see   \cite{cf:BBE,cf:BBD,cf:BDMbook} and references therein. 
We sketch a concise proof of Proposition~\ref{th:Y} because the context of \cite{cf:BBE,cf:BBD,cf:BDMbook} is really much wider and Proposition~\ref{th:Y} can be proven with substantially easier arguments. 

\medskip

We conclude this section by studying the asymptotic behavior of $F_\nu$.
From the results of
\cite{cf:BBD,cf:BDMbook} one can  extract that $F_\nu (x)=\nu((-\infty,x]) \sim  c x$ for $x \to \infty$;
however our techniques require a stronger result, see  \eqref{eq:Fnu-asympt} below, which requires substantially more constraining requirements than the hypotheses in \cite{cf:BBD,cf:BDMbook}.

Here is a preliminary result that allows the use of the Laplace transform.  Besides being an immediate consequence of the analysis in  \cite{cf:BBD,cf:BDMbook}, it can also be extracted from   \cite[Lemma~3.3]{cf:GGG}: since context and notations  here are slightly different, we give the proof in the appendix. We make a definite choice of $F_\nu$ by stipulating that 
$F_\nu\left(x_0\right) \, =\, 1$ for an $x_0$ which  is (strictly) inside the support of $\nu$.

\medskip

\begin{lemma}
\label{th:pre-bound}
There exists $c>0$ such that $F_\nu (x) \le c \exp(c x)$
for every $x\ge x_0$. 
\end{lemma}

\medskip

Here is the main result of this section:
\medskip

\begin{proposition}
\label{th:Fnu-asympt}
There exist three constants $m_\gz>0$, $c_\gz \in \bbR$  and $\varrho_\gz\in (0,1)$ such that for $x \to \infty$
\begin{equation}
\label{eq:Fnu-asympt}
F_\nu(x)\, =\, m_\gz x + c_ \gz + O\left( \exp(-\varrho_\gz x) \right)\, .
\end{equation}
Moreover, for every $\gd' \in (0 , \gd)$ (recall \eqref{eq:LR} for $\gd$) 
we have for $x \to -\infty$
\begin{equation}
\label{eq:Fnu-asympt2}
F_\nu(x)\, =\, O\left( \exp(\gd' x) \right)\, .
\end{equation}

\end{proposition}

\medskip

The constants carry the subscript $\gz$ because they are primarily determined by $\gz$, but $m_\gz$ and  $c_\gz$  depend also on the arbitrary choice of $x_0$:
$c_\gz/m_\gz$ and $\varrho_\gz$ instead depend only on $\gz$.

\medskip

\emph{Proof.} We start with \eqref{eq:Fnu-asympt}.
We are going to use the Laplace transform for (non negative) functions and measures supported on $[0, \infty)$:
$\widehat f(u) :=  \int_\bbR \exp(-ux) f(\dd x)=\int_{[0,\infty)} \exp(-ux) f(\dd x)$ and, if $f$ is absolutely continuous,
$\widehat f(u) = \int_0^\infty  \exp(-ux) f(x) \dd x$. 
If $F_f(x) \le \exp(c x)$ for $x$ large, then  $\widehat f(u)$ is analytic in the complex half  plane $ \Re u>c$ and $\vert \widehat f (u) \vert \le \widehat f  ( \Re u)$. 

In order to use the Laplace transform to study the right-tail of $\nu$ we introduce
$\tilde Y_n:= \log(1+ \exp(Y_n))$; then $(\tilde Y_n)$ is a Markov chain with
\begin{equation}
\label{eq:iterYtilde}
\tilde Y_{n+1}\, =\,  \log \left( 1+ \exp\left(  \logZ_{n+1} + \tilde Y_n\right)\right)\,=\,  \logZ_{n+1} + \tilde Y_n + \log \left( 1+ \exp\left( -\left( \logZ_{n+1} + \tilde Y_n\right)\right)\right)\, ,
\end{equation}
analogous to \eqref{eq:iterY}. 
The chains $(\tilde Y_n)$ and $(Y_n)$ are equivalent and the invariant measure $\mu$ of  $(\tilde Y_n)$ is directly related to $\nu$. In particular $\mu$ is $\gs$-finite,  $F_\mu(x)=0$ for every $x \le 0$ and,  
by Proposition~\ref{th:Y},  $F_\mu(x)< \infty$ for every $x$.  In fact, we have    
  $F_\nu(x)=F_\mu( \log(1+\exp(x)))$ so we can replace $F_\nu(x)$ with $F_\mu(x)$ in 
 \eqref{eq:Fnu-asympt} 
obtaining an equivalent expression.
Note the following characterization of $\mu$: for measurable $g: (0, \infty)\to [0, \infty)$
\begin{equation}
\label{eq:formu}
\int g \dd \mu \, =\, \int \int g\left( \log \left(1+ \exp(-z-x)\right) \right) \gz(\dd z) \mu(\dd x)\,.  
\end{equation}

By Lemma~\ref{th:pre-bound} we have that $\widehat F_\mu (u)$ is analytic in the domain $\Re u>c>0$.
We recall also that, in the same domain,  $\widehat \mu (u)=u\widehat F_\mu (u)$.
We call $b_\mu$ the \emph{optimal value of } $c$:
\begin{equation}
b_\mu:=\inf\{u\in \bbR :\,  \widehat \mu (u)< \infty\}\, ,
\end{equation}
so $b_\mu \le c$.
On the other hand Proposition~\ref{th:Y} tells us that $b_\mu \ge  0$, because $\nu(\bbR)=\mu([0, \infty))= \infty$. 
What we are going to show now is 
 that $b_\mu=0$.
 
 By 
\eqref{eq:formu} we readily see that 
\begin{equation}
\label{eq:fixed-Laplace}
\widehat \mu (u)\, =\, 
 \int \int (1+\exp(x+y))^{-u} \gz(\dd x) \mu ( \dd y)\,.
\end{equation}
We now use the identity \cite[(5.13.1)]{cf:DLMF}  (with $a=0$, $b=u$ and $s=w$, $a$, $b$ and $s$ are the notations in \cite{cf:DLMF})
\begin{equation}
\label{eq:Mellin-Barnes}
(1+z)^{-u}\, =\, \frac 1{ 2\pi i \gG(u)} \int_{w_0-i\infty}^{w_0 + i \infty} \gG (w) \gG(u-w) z^{-w} \dd w\,,
\end{equation}
which holds if  $w_0\in (0, \Re u)$ and  for every $z\in \bbC \setminus i(-\infty,0)$. We recall also the \emph{generalized  Stirling} formula \cite[(5.11.9)]{cf:DLMF}:
\begin{equation} 
\label{eq:5119}
\vert \gG(x+iy)\vert \stackrel {\vert y \vert \to \infty}\sim  \sqrt{2 \pi} \vert y \vert ^{x-1/2} \exp(-\pi \vert y \vert /2)\, ,
\end{equation}
for $\vert y \vert \to \infty$, uniformly for $x$ in bounded intervals.
Therefore by 
inserting \eqref{eq:Mellin-Barnes} into \eqref{eq:fixed-Laplace} and by using 
the  Fubini-Tonelli Theorem we obtain that for $b_\mu< w_0 < \Re u$  we have
\begin{equation}
\widehat \mu (u)\, =\, \frac{ 1}{2\pi i \gG(u)} 
 \int_{w_0-i\infty}^{w_0 + i \infty} \gG (w) \gG(u-w) \widehat \gz(w) \widehat \mu (w) \dd w\,.
\end{equation}
The next step is moving $w_0$ to the right of $\Re u$: more precisely to $(\Re u,\Re u+1)$. Given the decay for large imaginary part of the Gamma function, this can be done, but the stripe $\{z\in \bbC:\, w_0<\Re z <w_0+1\}$ contains
the simple pole of $\Gamma(w-u)$ at $w=u$. Since the residue of the integrand (times the prefactor) is 
$-\widehat \gz(u) \widehat \mu (u)$ and by remarking that we are performing the integral in the clockwise sense,  we obtain that for $w_0\in (\Re u,\Re u +1)$ 
\begin{equation}
\widehat \mu (u)\, =\, \widehat \gz(u) \widehat \mu (u) + \frac{ 1}{2\pi i \gG(u)} 
 \int_{w_0-i\infty}^{w_0 + i \infty} \gG (w) \gG(u-w) \widehat \gz(w)\widehat \mu (w) \dd w\,, 
\end{equation}
which, bringing the first term to the left hand side and dividing by $1-\widehat \zeta(w)$, we rewrite as
\begin{equation}
\label{eq:nuu}
\widehat \mu (u)\, =\, \frac{1}{2\pi i \gG(u)(1-\widehat \gz(u))} 
 \int_{w_0-i\infty}^{w_0 + i \infty} \gG (w) \gG(u-w)  \widehat \gz(w)\widehat \mu (w) \dd w\, =:\, P(u) I(u)\, ,
\end{equation}
where $P(u)$ denotes the pre-factor and $I(u)$ the integral in the second expression. 
We remark that this form makes it possible to analytically continue $I (u)$ to smaller values of $\Re u$, 
even down to $\Re u < w_0-1$ and $w_0$ is chosen larger (but arbitrarily close to $b_\mu$). 
But $P(u)$ has a singularity at  $u=0$ as we are going to explain:
\smallskip

\begin{itemize}
\item
$ \gG(u)\sim 1/u$ has a simple pole at $u=0$ (like for all the negative integers) and these are the only singularities; 
moreover the $\Gamma$ function has no zero.
\item
  $u \mapsto 1-\widehat \gz(u)$ is an entire function with
$1-\widehat \gz(u)\sim - \gs^2 u^2/2$,  $\gs^2:= \int x^2 \gz(\dd x)>0$: so, $1-\widehat \gz(u)$ has a double zero at the origin. We remark also that (with standard probabilistic notation)
$\widehat \gz(-it)=\gp_\logZ(t)$, i.e.\ it  is the characteristic function of a continuous random variable which yields   $\vert \widehat \gz(it)\vert< 1$
for $t\neq 0$ ($\vert\gp_\logZ(t)\vert=1$ for a $t\neq 0$ directly implies that $\logZ$ is discrete). By the Riemann-Lebesgue Lemma  we have $\vert \widehat \gz(it)\vert=o(1)$  for $\vert t \vert $
large. Therefore there exists $\varrho>0$ such that $1-\widehat \gz(u)=0$ for $\Re u >- \varrho$ only for $u=0$. We assume from now on that $\varrho \le 1$ (in any case,  $I(\cdot)$  has  been analytically continued only up to $\Re u >- 1$).
\end{itemize}

\smallskip

So $P(\cdot)$ does contribute a simple pole at $z=0$: a priori there is still the possibility that $I(0)=0$ making this singularity removable, but in fact this is not possible because we have already remarked that $b_\mu \ge 0$. 

\smallskip

We have therefore proven not only that $b_\mu=0$, but also that
$\widehat \mu $ can be meromorphically extended to $\{z \in \bbC:\, \Re z >-\varrho \}$ and the only singularity in this region is the simple pole at zero
\begin{equation}
 \widehat \mu (u)   \,=\,  \frac {m_\gz} u +c_\gz + u H(u) \, ,
\end{equation}
with $m_\gz>0$ because $\widehat\mu(u)>0$ for $u>0$, $c_\gz \in \bbR$ and $H(\cdot)$ is an analytic function on the domain $\{z \in \bbC:\, \Re z >-\varrho \}$. This of course tells us that also 
$\widehat F_\mu$ can be meromorphically extended to the same domain and 
\begin{equation}
 \widehat F_\mu(u)   \,=\,  \frac {m_\gz} {u^2} +\frac{c_\gz} u +  H(u) \, .
\end{equation}
\medskip

 
 We now use the classical Mellin-Bromwich formula for the inverse Laplace transform: for every $b>b_\mu=0$ 
 \begin{equation}
 \label{eq:MeBr}
 F_\mu(x)\, =\,\frac {1}{ 2\pi i}\lim_{R \to \infty} \int_{b-i R}^{b+ iR} \widehat F_\mu(u) \exp(ux) \dd u\, .
 \end{equation}
 We introduce the rectangular contour $\mathtt{C}_R$ made of the segment of integration 
 in \eqref{eq:MeBr}, the segment $\{u: \, \Re u =-\varrho_\gz$ and $\Im u \in[-R, R]\}$ and by the two segments $\{u:\, \Im u = \pm R$ and $\Re u\in [-\varrho_\gz, a]\}$. The orientation is counter-clockwise and $\varrho_\gz $ is chosen in $(0, \varrho)$.
 The integration along this contour is equal to the residue of the double pole at the origin:
\begin{equation}
\widehat F_\mu(u) \exp(ux)\, =\,  \left(\frac {m_\gz} {u^2} +\frac{c_\gz} u +  O(1)\right)( 1+ xu + O(u^2))\,=\,
\frac {m_\gz} {u^2}
+ \frac {c_\gz +m_\gz x} u + O(1)\, ,
\end{equation}
hence the residue is $c_\gz +m_\gz x$ and 
\begin{equation}
\frac {1}{ 2\pi i} \int_{\mathtt{C}_R} \widehat F_\mu(u) \exp(ux) \dd u\, =\, c_\gz +m_\gz x\, .
\end{equation}

We are therefore left with  showing that:

\smallskip

\begin{enumerate}
\item
 the contribution of the integrals along the segments  
 $\{u:\, \Im u = \pm R$ and $\Re u\in [-\varrho_\gz, a]\}$ vanishes as $R\to \infty$;
 \item the contribution of the integrals along $\{u: \, \Re u =-\varrho_\gz$ and $\Im u \in[-R, R]\}$ is $O(\exp(-\varrho_\gz x))$ for $x \to \infty$. 
\end{enumerate}

 \smallskip 
 
 Both estimates rely on 
  
\medskip

\begin{lemma}
\label{th:gGctrl}
We consider $w_0>0$  and a closed interval $I\subset(w_0-1/2, w_0)$. For $ \vert w_1 \vert\to \infty$ and uniformly in $\Re u\in I$
\begin{equation}
\label{eq:gGctrl}
 \frac{1}{ \left \vert\gG(u)\right\vert} 
 \int_{w_0-i\infty}^{w_0 + i \infty} \left \vert  \gG (w)\right \vert
  \left \vert \gG(u-w) \right \vert \dd w\, =\, O \left( \left \vert \Im u \right \vert^{1/2}\right)
\end{equation}
\end{lemma}

\medskip

Lemma~\ref{th:gGctrl} is \cite[Lemma~4.4]{cf:GGG}, but the proof can also be obtained directly from \eqref{eq:5119}. We remark that the dominant contribution comes from integrating for $w$ in a neighborhood of the real axis.

\smallskip

Moreover for both (1) and (2) the crucial formula is \eqref{eq:nuu} because, in conjunction with \eqref{eq:gGctrl} and 
recalling that $\widehat F_\mu(u)= \widehat \mu (u) /u$,  it tells us 
that 
\begin{equation}
\vert \widehat F_\mu(u) \vert \,=\, O\left( \big\vert \widehat \gz(u)\big\vert \left \vert \Im u \right \vert^{-1/2}
\right)\, ,
\end{equation}
for $\vert \Im u \vert$ large and uniformly for $\Re u$ in the interval we consider. 

Therefore the contour integration in point (1) is $O(1/R^{1/2})$ and for what concerns point (2) we have that for every $R$
\begin{equation}
\left \vert 
 \int_{-\varrho_\gz -iR}^{-\varrho_\gz +iR} \widehat F_\mu(u) \exp(ux) \dd u\right\vert \, \le \, 
\exp(- \varrho_\gz x) \int_{-\infty}^\infty \left \vert \widehat F_\mu\left(-\varrho_\gz +i y \right) \right\vert \dd y\, , 
\end{equation} 
with the integral in the right-hand side that converges because of \eqref{eq:Fourier1/2}. This completes the proof of \eqref{eq:Fnu-asympt}.
 
 \smallskip
 
 For \eqref{eq:Fnu-asympt2} we use \eqref{eq:fornu} 
 with $g=\ind_{(-\infty, x]}$, that is
 \begin{equation}
 \label{eq:Fnu-c}
 F_\nu(x)\,=\, \int \int \ind_{(-\infty, x]} (z+h(y)) \nu(\dd y) \gz (\dd z)\, .
 \end{equation}
 By splitting the integral in the $y$ variable into positive and negative values and by using that $h(y) \ge 0$ for $y$ negative and $h(y)\ge y$ for $y$ positive we see that
\begin{equation}
\begin{split}
 F_\nu(x)\,&\le \, F_\nu(0) F_\gz(x)+ \int \int \ind_{(-\infty, x]} (z+y) \ind_{(0, \infty)}(y) \nu(\dd y) \gz (\dd z)
 \\
 &\le \,  F_\nu(0) F_\gz(x)+  \int \int \ind_{(-\infty, x-z]} (y) \ind_{(0, \infty)}(x-z) \nu(\dd y) \gz (\dd z)
 \\
 = \,  F_\nu(0) &F_\gz(x)+   \int F_\nu (x-z) \ind_{(0, \infty)}(x-z)  \gz (\dd z) 
 \le F_\nu(0) F_\gz(x)+ C\int_0^\infty y \gz(x-y) \dd y
 \, ,
 \end{split}
 \end{equation} 
where in  the last step $C$ is a positive constant and we have used \eqref{eq:Fnu-asympt}.
But $\int_0^\infty y \gz(x-y) \dd y= \int_{-\infty}^x F_\gz (z) \dd z$ and \eqref{eq:Fnu-asympt2}  is established.

The proof of Proposition~\ref{th:Fnu-asympt} is therefore complete.
\qed

\medskip

\begin{rem}
\label{rem:nu}
 \eqref{eq:Fnu-c} also characterizes $\nu_k$ if we replace $h(\cdot)$ with $h_k(\cdot)$. From this we also have the characterization
\begin{equation}
\label{eq:remnu}
F_{\nu_k}(x)\, =\, \int_ \bbR F_\gz(x-h_k(y)) \nu_k (\dd y)\, ,
\end{equation}
from which we see that $F_{\nu_k}$ is $C^1$, so $\nu_k$ has a density,  because $0 \le \gz(x-h_k(y)) \le \Vert \gz\Vert_\infty$. 

Similarly $\nu$ has a density, but the argument needs to be refined.
Equation~\eqref{eq:remnu} still holds also without the subscripts $k$,
but because $\nu(\bbR)=\infty$
we need an appropriate upper bound for $ \gz(x-h(y))$ in order to take the derivative under the integral. 
Since $\nu((-\infty,x])< \infty$ for every $x$, $\gz \in \bbL^\infty$ reduces the problem to finding and upper bound for $y>0$ and large (of course we can focus on $x$ in a compact set). But this follows because there exists
$\gep>0$ and $z_0$ such that $\gz(z) \le \exp( \gep z)$ for $z \le z_0$. In order to show this 
and  to make the argument  more readable, let us replace $\gz(z)$ with $\gz(-z)$. If $\gz(z) \le \exp( -\gep z)$ for $z$
larger than some value $z_0$ is false, then there exists a sequence  $(z_j)$ with $\lim_j z_j=\infty$ and $\gz(z_j)> \exp(-\gep z_j)$.
Uniform $\theta$-H\"older continuity implies that $\gz(z) \ge  \exp(-\gep z_j)/2$ for $\vert z-z_j\vert \le ((e^{-\gep z_j})/2C)^{1/ \theta}$. So, possibly for $j$ large,  $\int_{z_j/2}^\infty \gz(z) \dd z \ge  \exp(-\gep z_j) ((e^{-\gep z_j})/2C)^{1/ \theta}$.
Since the left-hand side decays exponentially, by choosing $\gep>0$ small we reach a contradiction.
So we have shown that  $\gz$ vanishes exponentially fast: using that  $\nu((-\infty, x])$ has linear asymptotic growth, $F_\nu\in C^1$ follows.
\end{rem}

\section{The DH probability, the one step bound   and the proof the main result}
\label{sec:DH}

We now  define the DH probability, i.e.\ the measure that we expect (and  will prove)  to be close to  the invariant probability $\nu_k$. 
It is built by gluing together the invariant measure for the edge process at $k$ and for the one at $-k$.  
Unless $\gz$ is symmetric, the two edge limit problems are not the same, since the one on the left involves 
$x \mapsto\gz(x)$ as jump density probability and the one one the right involves $x \mapsto\gz(-x)$. So the cumulative function on the left (respectively, right) edge will be denoted by $F_\lar(\cdot)$ (respectively,
$F_\rar(\cdot)$). Moreover we choose to normalize the cumulative functions so that they are, for $x \to \infty$, equivalent to $x$.
Proposition~\ref{th:Fnu-asympt}  then says that there exist $c_s\in \bbR$, $s$ is $\lar$ or $\rar$, and  $\varrho\in (0,1)$ such that for $x \to \infty$
\begin{equation}
\label{eq:asymptlrar}
F_s(x)= x+ c_s + O\left(\exp(-\varrho x)\right)\, .
\end{equation}
We define the DH probability by giving its 
  integrated tail probability $G_k(x)= \nu^{\textrm{DH}}_k((x, \infty))$ (cf. Sec.~\ref{sec:strategy}):
\begin{equation}
G_k(x)\, :=\, \begin{cases} 
F_{\rar}(k-x) / C_k  & \text{ if } x \ge 0\, ,
\\
1- \left(F_\lar(x+k)/C_k\right) & \text{ if } x \le 0\, ,
\end{cases}
\end{equation}
where the fact that the definition  must be consistent at $x=0$ 
fixes the value of $C_k>0$ which, therefore, by \eqref{eq:asymptlrar}, satisfies
\begin{equation}
\label{eq:Ck}
C_k\, =\, 2k +c_\lar +c_\rar + O(\exp(-\varrho k))\,.
\end{equation}  
We register also that for the cumulative probability  $F_k(\cdot)=1-G_k(\cdot)$ we have
\begin{equation}
\label{eq:cumFk}
F_k(x)\, :=\, \begin{cases} 
1-\left(F_{\rar}(k-x) / C_k\right)  & \text{ if } x \ge 0\, ,
\\
F_\lar(x+k)/C_k & \text{ if } x \le 0\, .
\end{cases}
\end{equation}
The next fact is straightforward  computation, but it is of course central   for us:

\medskip

\begin{lemma}
\label{th:Lyap}
For $k \to \infty$
\begin{equation}
\label{eq:Lyap}
L_k[G_k]\,=\, \frac 1{C_k}\int_\bbR \frac{F_\rar(y)}{1+e^y} \dd y+ O( \exp(-k))\,=\,   \frac 1{C_k}\int_\bbR \frac{F_\lar(y)}{1+e^y} \dd y+ O( \exp(-k))\,.
\end{equation}
\end{lemma}

\medskip

\noindent
\emph{Proof.} By the first equality in \eqref{eq:L_k} we have
\begin{equation}
\label{eq:for-Lyap-1}
L_k[G_k]\,=\, \frac 1{C_k}\int_0^\infty \frac{F_\rar(k-x)}{1+e^{k-x}} \dd x+
\int_{-\infty}^0 \frac{G_k(x)} {1+e^{k-x}} \dd x\, .
\end{equation}
Using $G_k(x)\in [0, 1]$ we readily see that the second addendum is smaller than $2\exp(-k)$.
Moreover
\begin{equation}
\int_0^\infty \frac{F_\rar(k-x)}{1+e^{k-x}} \dd x\, =\, 
\int_{\bbR} \frac{F_\rar(y)}{1+e^{y}} \dd y - \int_k^\infty \frac{F_\rar(y)}{1+e^{y}} \dd y\, =\, 
\int_{\bbR} \frac{F_\rar(y)}{1+e^{y}} \dd y +O\left(k \exp(-k)\right)
\, ,
\end{equation}
where in the last step we used \eqref{eq:asymptlrar}. 

Repeating the same argument starting with the second equality in \eqref{eq:L_k} gives the same expression with $F_\lar$ replacing $F_\rar$.
\qed

\medskip

\begin{rem}
\label{rem:lyapsim}
\eqref{eq:Lyap} of course implies $\int (1+e^y)^{-1}  F_s(y) \dd y $ does not depend on $s$. Equivalently, by integration by parts:
\begin{equation}
\label{eq:Lyap=}
\int_\bbR \log\left( 1+e^{-y} \right) \nu _\rar (\dd y) \, =\, 
\int_\bbR \log\left( 1+e^{-y} \right) \nu _\lar (\dd y)\,. 
\end{equation}
One way to see this directly is to observe that by \eqref{eq:XMC} and \eqref{eq:map} 
we have $\int_\bbR (x-h_k(x)) \nu_k(\dd x)=0$, with $\nu_k$ the invariant probability of the Markov chain defined by \eqref{eq:XMC}. But this is equivalent to 
\begin{equation}
\int_\bbR \log \left( 1+ e^{k-x} \right) \nu_k(\dd x)\, =\, \int_\bbR \log \left( 1+ e^{-k-x} \right) \nu_k(\dd x)\, ,
\end{equation}
and \eqref{eq:Lyap=} follows by exploiting 
the  strong form of convergence -- see Remark~\ref{rem:vague} -- we have for $2k \nu_k \Theta_{-k}^{-1}$ towards $\nu_\lar$ and for  $2k \nu_k \Theta_k^{-1}$ towards $\nu_\rar$: here $\Theta_a(x):= x+a$. 
\end{rem}


\medskip

\begin{proposition}
\label{th:onestep}
There exists $\eta>0$ such that
\begin{equation}
\label{eq:2addenda}
\left\Vert T G_k - G_k \right \Vert_1\, =\, O\left( \exp(-\eta k) \right)\,.
\end{equation}
\end{proposition}

\medskip

The proof is written so that the constant $\eta$ in \eqref{eq:2addenda} can be chosen arbitrarily in $(0, \min(\gd/2, \theta)/2, \varrho)$,
with $\gd$ in  \eqref{eq:LR}, $\theta$ in  \eqref{eq:Holder} and $\varrho$ is the constant $\varrho_\zeta$ in Proposition~\ref{th:Fnu-asympt}: a slight modification of the proof yields that we can choose $\eta \in (0, \min(\gd, \theta, \varrho))$. 

\smallskip 

Before proving Proposition~\ref{th:onestep} \emph{we remark in an official way} that it is the last brick needed for our main result. 

\medskip

\noindent
\emph{Proof of Theorem~\ref{th:main}}. Recall that $\nu_k$ is the invariant probability of the main chain $X$.
It suffices to 
 choose $G_\gamma=G_k$ and apply Corollary~\ref{th:bound-cor}. By using Proposition~\ref{th:onestep}
we obtain that 
\begin{equation} 
\label{eq:bound-cor2}
\left \vert \cL(k)- L_k[G_k] \right \vert \, \le\, \left \Vert G_{\nu_k} - G _ k \right \Vert _1\, = \,O\left( k^2 \left( \log k \right)^C \exp(-\eta k) \right)\, .
\end{equation}
We then conclude by Lemma~\ref{th:Lyap} and by exploiting the expression  \eqref{eq:Ck} for $C_k$. We
obtain
\begin{equation}
\label{eq:constants}
\kappa_1\, =\, \frac 12 \int_\bbR \frac{F_{\lar}(y)}{1+e^y} \dd y\, =\, \frac 12 \int_\bbR \frac{F_{\rar}(y)}{1+e^y} \dd y \ \ \text{ and } \ \ \kappa_2 = \frac 12 \left(c_\lar +c_\rar\right)\, .
\end{equation}
\qed
\medskip

\begin{rem}
\label{rem:vague}
We remark that a byproduct of the $\bbL^1$ control in \eqref{eq:bound-cor2} is that 
$2k \nu_k \theta_{-k}^{-1}$ converges vaguely for $k \to \infty$ towards $\nu_\lar$ and  $2k \nu_k \theta_k^{-1}$ converges vaguely towards $\nu_\rar$. Vague convergence, i.e.\ with test functions that are compactly supported and $C^0$,  is an elegant statement, but \eqref{eq:bound-cor2}  is much stronger and, in particular, it solves the issue raised in Remark~\ref{rem:lyapsim}.
\end{rem}

\medskip

\noindent
\emph{Proof of Proposition~\ref{th:onestep}.}
Of course
\begin{equation}
\label{eq:2addenda-2}
\left\Vert T G_k - G_k \right \Vert_1\, =\, 
\left\Vert (T G_k - G_k )\ind_{(-\infty,0)}\right \Vert_1 + \left\Vert (T G_k - G_k )\ind_{(0, \infty)}\right \Vert_1\, ,
\end{equation}
and, even if the two terms on the right may be (and typically are) different unless $\gz$ is symmetric, 
they can be treated in the same way because
both of them can be written as 
\begin{equation}
A_{k,1}\, :=\,\int_{-\infty}^0 \left \vert  \int_ \bbR F(y) h'_k(y) \gz( x-h_k(y)) \dd y - F(x) \right\vert  \dd x\, , 
\end{equation}
 with $F(\cdot)$ which is $F_k(\cdot)$, given in \eqref{eq:cumFk},  for the first 
 addendum in  the right-hand side of \eqref{eq:2addenda-2}, and $F(\cdot)$ that is instead replaced by the right-hand side in \eqref{eq:cumFk} with $ \lar$ and $\rar$ exchanged.
 Arbitrarily, we choose to work with the first addendum and the bound we are after is achieved in two steps.

The first step in controlling $A_{k,1}$ is to remark that we can avoid the nuisance of the fact that 
$F_k(y)$ has two different expressions according to the sign of $y$. Namely, we want to switch to:
\begin{equation}
A_{k,2}\, :=\,\frac 1{C_k}\int_{-\infty}^0 \left \vert  \int_ \bbR F_\lar(y+k) h'_k(y) \gz( x-h_k(y)) \dd y - F_\lar(x+k) \right\vert  \dd x\, .
\end{equation}
We can do this because 
\begin{equation}
\left\vert A_{k,2}-A_{k,1}\right\vert\, \le \,
\int_{-\infty}^0    \int_{0}^{\infty} \left \vert \frac{
F_\lar(k+y) + F_\rar(k-y)}{C_k}-1 \right \vert  h'_k(y) \gz( x-h_k(y)) \dd y \dd x\, ,
\end{equation}
and, by \eqref{eq:asymptlrar} and \eqref{eq:Ck}, for $y\ge 0$ we have
\begin{equation}
\left \vert \frac{F_\lar(k+y) + F_\rar(k-y)}{C_k}-1 \right \vert \,\le\, \frac{C}{C_k} \times  \begin{cases}\exp(-\varrho (k-y)) & \text{ if } y\in [0, k] \, ,
\\
y-k+1 & \text{ if } y >k\, ,
\end{cases}
\end{equation}
for a suitably chosen $C>0$. So 
 $\vert A_{k,2}-A_{k,1}\vert$ is bounded above by 
\begin{equation}
\label{eq:interm2.k}
     \int_{0}^{k}  \exp(-\varrho (k-y)) F_\gz(-h_k(y)) h'_k(y)\dd y +
   \int_{k}^{\infty} (y-k+1)  F_\gz(-h_k(y)) h'_k(y) \dd y\,,
 \end{equation}
 times $C/C_k$.
 So, by keeping in mind that $F_\zeta(-x)=O(\exp(-\gd x)$, that $h'_k(y) =O(\exp(-(y-k)))$ for $y-k \to \infty$, that  $ h'_k(y)\in (0,1)$ and by remarking that  $y-h_k(y)\le \log 2$ for $y\in [0,k]$,  with adequate choice of $C$
 \begin{equation}
  \left\vert A_{k,2}-A_{k,1}\right\vert\, \le  \,  \frac{C}{C_k} \left( 
   \int_{0}^{k}  e^{-\varrho (k-y)} e^{-\gd y} \dd y + e^{-\gd k}
   \int_{k}^{\infty} (y-k+1)   e^{-(y-k)} \dd y\right)\, =\, O\left( {e^{-(\gd \wedge \varrho) k}}\right),
  \end{equation}
  with $\gd \wedge \varrho= \min(\gd, \varrho)$.

The second step is the control of  $A_{k,2}$. For this we 
introduce
\begin{equation}
\mathtt{h}_k(x)\, :=\, h(x+k)-k\, 
\end{equation}
 and remark that
 \begin{equation}
\label{eq:ident-before-asympt1}
0 \, \le \, \mathtt{h}_k(x)- h_k(x)\, =\, \log(1+ \exp(x-k))\,= \begin{cases}
O(\exp(x-k))& \text{ for } x-k\to -\infty\, , \\
O(x-k) & \text{ for } x-k\to +\infty\, ,
\end{cases}
\end{equation}
\begin{equation}
\label{eq:ident-before-asympt2}
0 \, \le \, \mathtt{h}'_k(x)- h'_k(x)\, =\, 1/(1+ \exp(k-x))\, = \begin{cases}
O(\exp(x-k))& \text{ for } x-k\to -\infty\, , \\
\le 1 & \text{ for every } x \text{ and } k\, .
\end{cases}
\end{equation}
 Moreover we stress that $h_k'(x)$ and $\mathtt{h}'_k(x)$ are in $(0,1)$ for every $x$ and, with reference to Figure~\ref{fig:1}, the two functions almost coincide for $x\ll k$. They start to differ when $x$ approaches $k$ from the left because  $\mathtt{h}'_k(x)$ keeps being very close to one and $\mathtt{h}_k(x)\sim x$. 
 
 Since by the characterizing property of $F_\lar(\cdot)$
 \begin{equation}
 \label{eq:passo-interm1}
 \int_ \bbR F_\lar(y+k) \mathtt{h}'_k(y) \gz( x-\mathtt{h}_k(y)) \dd y - F_\lar(x+k)=
  \int_ \bbR F_\lar(y)h'(y) \gz( x-h(y)) \dd y - F_\lar(x)\, =\,0 ,
 \end{equation}
 we have 
 \begin{equation}
 A_{k,2}\, =\, \frac 1{C_k} \int_{-\infty}^0  \left\vert 
   \int_ \bbR F_\lar(y+k) \left(
   h'_k(y) \gz( x-h_k(y))-   \mathtt{h}'_k(y) \gz( x-\mathtt{h}_k(y)) 
\right)    \dd y  \right\vert \dd x\, ,
 \end{equation}
 We split the integral over $\bbR$ according to $ y < -3k/2$, $y\in[-3k/2, k/2]$ and $y>k/2$. We have 
 \begin{multline}
  \int_{-\infty}^0 \left \vert
  \int_{k/2}^\infty F_\lar(y+k) \left(
   h'_k(y) \gz( x-h_k(y))-   \mathtt{h}'_k(y) \gz( x-\mathtt{h}_k(y)) 
\right)    \dd y  \right\vert \dd x
\, \le \\
C
 \int_{-\infty}^0 
  \int_{k/2}^\infty (y+k) \left(
   h'_k(y) \gz( x-h_k(y)) +  \gz( x-\mathtt{h}_k(y)) 
\right)    \dd y  \dd x\, =
\\
C
  \int_{k/2}^\infty (y+k) \left(
   h'_k(y) F_\gz( -h_k(y)) +  F_\gz( -\mathtt{h}_k(y)) 
\right)    \dd y  \, \le \, 
C'
  \int_{k/2}^\infty y \left(
   e^{(k-y)\wedge 0 -\frac\gd 2 k}  +  e^{-\gd y}
\right)    \dd y ,
 \end{multline}
 which is $O(k^2 \exp(-\gd k/2))$. Similarly with $\gd'\in (0, \gd)$ (recall \eqref{eq:Fnu-asympt2})
 \begin{multline}
  \int_{-\infty}^0 \left \vert
  \int_{-\infty}^{-3k/2} F_\lar(y+k) \left(
   h'_k(y) \gz( x-h_k(y))-   \mathtt{h}'_k(y) \gz( x-\mathtt{h}_k(y)) 
\right)    \dd y  \right\vert \dd x
\, \le \\
C  \int_{-\infty}^{-3k/2} e^{\gd' (y+k)}e^{y+k} (F_\gz(-h_k(y))+F_\gz( -\mathtt{h}_k(y))\dd y\, \le \, 
2C  \int_{-\infty}^{-3k/2} e^{\gd' (y+k)}e^{y+k} \dd y\,
 \,,
\end{multline}
which is $O(\exp(-k(1+ \gd')/2))$.
We have therefore obtained that $A_{k,2}$ is equal to 
\begin{equation}
  \frac 1{C_k} \int_{-\infty}^0  \left\vert 
   \int_{-3k/2}^{k/2} F_\lar(y+k) \left(
   h'_k(y) \gz( x-h_k(y))-   \mathtt{h}'_k(y) \gz( x-\mathtt{h}_k(y)) 
\right)    \dd y  \right\vert \dd x  + O\left(k e^{- \frac {\gd'}  2 k}\right)\, .
 \end{equation}
For the last estimate we use 
 \begin{multline}
  \label{eq:passo-interm2}
  \left\vert 
    \int_{-3k/2}^{k/2} F_\lar(y+k) \left(
   h'_k(y) \gz( x-h_k(y))-   \mathtt{h}'_k(y) \gz( x-\mathtt{h}_k(y)) 
\right)    \dd y  \right\vert \, \le 
 \\
  \int_{-3k/2}^{k/2}  F_\lar(y+k)\big(\underbrace{\left\vert 
  h'_k(y) -   \mathtt{h}'_k(y)\right \vert \gz( x-h_k(y))}_{1} +   \underbrace{\mathtt{h}'_k(y)  \left\vert \gz( x-h_k(y))-  \gz( x-\mathtt{h}_k(y))  
\right\vert}_2 \big)   \dd y\, ,
 \end{multline}
 and we write the last line as  $T_1(x)+T_2(x)$ in the obvious way. Since $F_\lar(y+k) \le 2k$ for $y\le k/2$
 \begin{equation}
 \begin{split}
 \int_{-\infty}^0 T_1(x) \dd x \, &=\,   \int_{-3k/2}^{k/2}  F_\lar(y+k)\left\vert 
  h'_k(y) -   \mathtt{h}'_k(y)\right \vert F_\gz( -h_k(y)) \dd y\\
  &\le \, 2 k
   \int_{-3k/2}^{k/2}  \left\vert 
  h'_k(y) -   \mathtt{h}'_k(y)\right \vert  \dd y 
  \, \le \,  C k
   \int_{-3k/2}^{k/2}  e^{y-k}  \dd y \, \le \,  O\left( k \exp(- k/2)\right)\,.
 \end{split}
 \end{equation}
 Moreover
 \begin{equation}
 \begin{split}
 \int_{-\infty}^0 T_2(x) \dd x \, &=\,    \int_{-\infty}^0\int_{-3k/2}^{k/2}  F_\lar(y+k)  \mathtt{h}'_k(y)  \left\vert \gz( x-h_k(y))-  \gz( x-\mathtt{h}_k(y))  
\right\vert \dd y \dd x 
\\
& \le \, 
 \int_{-2k}^0\int_{-3k/2}^{k/2}  F_\lar(y+k)  
 C_\gz \left\vert h_k(y))-  \mathtt{h}_k(y)
\right\vert^\theta \dd y \dd x\\
& \ \ \ + C \int_{-3k/2}^{k/2}  F_\lar(y+k)  
\left(F_\gz( -2k-h_k(y))+ F_\gz( -2k-\mathtt{h}_k(y))  
\right) 
 \dd y
 \\
  \le  C' k \int_{-2k}^0  &\dd x \int_{-3k/2}^{k/2}   
 \exp(\theta(y-k)) \dd y + 8Ck^2 F_\gz( -k)\, =\, O\left( k^2 \exp\left(-\left(\frac \theta 2\wedge \gd\right) k\right) \right),
\end{split}
 \end{equation}
 where in the second step we have exploited the regularity of $\gz$ for $x\in (-2k, 0)$ --
 the constant $C_\gz$ is the left-hand side in \eqref{eq:Holder} -- and, for smaller $x$'s we have just bounded $ \vert \gz( x-h_k(y))-  \gz( x-\mathtt{h}_k(y))\vert$ with
 $  \gz( x-h_k(y))+  \gz( x-\mathtt{h}_k(y))$ and we have performed the integral over $x$, using once again $ F_\lar(y+k)\le 2k$ in the range we consider.

The proof of Proposition~\ref{th:onestep} is therefore complete.
\qed

\section{Diffusion estimate: the proof of Proposition~\ref{th:bound}}
\label{sec:bound}

In this section we essentially estimate the speed of convergence to stationarity of the chain $X$. The process is essentially a symmetric random walk far from the boundary, but the positive recurrent character crucially depends on the visits to the boundary. By a preliminary manipulation we reduce the estimate we need to estimating the expectation of the time of first exit from $(-k,k)$ and this is achieved by  diffusion approximation 
on $(-k+ 2 \log \log k, k-2 \log \log k)$ and by a rough estimate on the  remaining portions of length $2\log \log k$. 

\medskip

\noindent
\emph{Proof of Proposition~\ref{th:bound}.} 
In order to work with more explicit constants we give the proof under the assumption that 
both $\bbP(\logZ > \log 2)>0$ and $\bbP(\logZ < -\log 2 )>0$: we explain in Remark~\ref{rem:2} how the proof can be easily generalized.

\medskip
Our proof involves two auxiliary Markov chains, closely related to $X$, which we now introduce.

The first one is  $X^\sim=(X^\sim_n)_{n\in\bbN_0}$ with state space $\bbR\cup \{\cC\}$. For $X^\sim$  the state $\cC$ is absorbing (\emph{cemetery}), that is $X^\sim_n=\cC$ implies $X^\sim_m=\cC$ for every $m$ larger than $n$. 
On the other hand if $X^\sim_n =y \in \bbR$ the probability that $X^\sim_{n+1}= \cC$ is $1-h'_k(y)$.
Therefore if $X^\sim_n =y \in \bbR$ the chain \emph{survives} with probability $h'_k(y)$ and, in this case, it jumps to $x$ 
with transition density $\gz(x-h_k(y))$. Note that $h'_k(x)$ is even, decreasing for $x>0$,
$h'_k(x)= O(\exp(-k+x))+ O(\exp(-k-x))$ and 
\begin{equation}
\label{eq:1/2}
h'_k(\pm k)\, =\, \frac 1 2 - \frac{1}{e^{2k}+1} \, < \, \frac 12\,,
\end{equation}
so $h'_k(x) < 1/2$ if $\vert x\vert \ge k$, and $h'_k(x) \le h'_k(0)=1$ for every $x$. 

\medskip

By the observations we just made it is easily seen that  the chain $X^\sim$ can be coupled with a  simpler chain $X^\approx=
( X^\approx_n)_{n\in\bbN_0}$ that has a smaller \emph{death probability}.
The transition kernel for $X^\approx$ from $y\in (-k,k)$ to $x\in \bbR$ is given by
$\gz(x-h_k(y))$ and it is therefore a probability: transitions from $(-k,k)$ to $\cC$ are forbidden.
On the other hand the kernel from $y \in \bbR \setminus (-k,k)$ to $x \in \bbR$ is $\gz(x-h_k(y))/2$ and  
the probability of going from $y \in \bbR \setminus (-k,k)$
to $\cC$ is $1/2$. $\cC$ is absorbing also for this chain.

\medskip

Of course these two chains are \emph{dominated} by the even simpler chain $X=
( X_n)_{n\in\bbN_0}$ for which the death probability is zero: its natural state space is just $\bbR$, but of course 
we can keep $\bbR \cup \{\cC\}$ as state space and the absorbing state does not communicate with $\bbR$,  and the transition kernel from $y$ to $x$ is the probability density $\gz(x-h_k(y))$. In fact $X$ is just the main Markov chain we consider in this walk, i.e.\ \eqref{eq:XMC}. 

We denote by $\bbP_y$ the (joint) law of the chains $X$, $X^\sim$ and $X^\approx$ and the index $y$
denotes the (common) initial condition.

\medskip

Let us assume  that $G(\cdot)$ is non negative (the general result is recovered by linearity).
Of course
\begin{equation}
\int_\bbR T_0G(x) \dd x \, =\, \int_{\bbR} h'_k(y) \gz(x-h_k(y)) G(y) \dd y\, ,
\end{equation}
and a  direct computation shows that for $n=2,3, \ldots$
\begin{equation}
\int_\bbR  T_0^n G(x) \dd x \, =\, \int_{\bbR^n} h_k'(y_1) 
\left(
\prod_{j=1}^{n-1} \gz(y_j-h_k(y_{j+1}))  h'_k (y_{j+1})\right)
\dd y_1 \ldots \dd y_n\, .
\end{equation}
Therefore for $n\in \bbN$ we have the probabilistic representation and bound:
\begin{equation}
\begin{split}
\int_\bbR  T_0^n G(x) \dd x \, =\,
\int_\bbR G(y)\bbP_y\left( X^{\sim}_n\not= \cC \right) \dd y \, &\le \, 
\int_\bbR G(y)\bbP_y\left( X^{\approx}_n\not= \cC \right) \dd y
\\
&=\, \int_\bbR G(y)\bbE_y\left[ 2^{-H(n)} \right] \dd y\, ,
\end{split}
\end{equation}
with $H(n):=\vert \{j=0,1, \ldots, n:\, \vert X_j\vert \ge k \}\vert$. Let us set $\tau^H_0=0$ and, for $j \in \bbN$,
$\tau^H_j:= \inf \{ n \in \bbN: \, H(n)=j\}$, so $\tau^H_j$ is the time of the $j^{\textrm{th}}$ entry of $X$
into $\bbR \setminus (-k,k)$. With these notations, by the Fubini-Tonelli Theorem, we have
\begin{equation}
\sum_{n =0}^\infty
\int_\bbR G(y)\bbE_y\left[ 2^{-H(n)} \right] \dd y \, =\, 
\int_\bbR G(y)
\bbE_y\left[  \sum_{n =0}^\infty2^{-H(n)} \right]  \dd y \, ,
\end{equation}
and, again  by the Fubini-Tonelli Theorem and using the definitions:
\begin{equation}
\bbE_y\left[  \sum_{n =0}^\infty2^{-H(n)} \right]\, =\, 
\sum_{j=0}^\infty 2^{-j}
\bbE_y\left[  \sum_{n =0}^\infty
\ind_{H(n)=j} \right]\, =\, \sum_{j=0}^\infty 2^{-j} \bbE_y\left[ \tau^H_{j+1}-  \tau^H_{j}\right]\, \le \, 
2 \sup_y  \bbE_y\left[ \tau^H_1\right]\,.
\end{equation}
For every $r>0$ we set $\tau_r:= \inf\{n\ge 0:\, X_n \in \bbR \setminus (-r,r)\}$.
So  $\tau^H_1=\tau_k$ and the steps we have just developed yield that for $G(\cdot)\ge 0$
\begin{equation}
\sum_{n=0}^\infty
\int_\bbR  T_0^n G(x) \dd x \, \le \, 2  \sup_y  \bbE_y\left[ \tau_k\right] \int_\bbR   G(x) \dd x \, ,
\end{equation}
which readily implies for $G \in L^1$
\begin{equation}
\label{eq:tbuwlemTk}
\sum_{n=0}^\infty
\left\Vert T_0^n G \right\Vert_1\, \le \, \sum_{n=0}^\infty
\int_\bbR  T_0^n \vert G(x)\vert \dd x \, \le \, 2  \sup_y  \bbE_y\left[ \tau_k\right] \Vert G \Vert_1 \, .
\end{equation}

We are going to show:

\medskip

\begin{lemma}
\label{th:Tk}
There exists $C>0$ and $k_0>0$ such that for every $k \ge k_0$
\begin{equation}
\label{eq:Tk}
\sup_{x }\bbE_x\left[ \tau_k\right] \, \le \, \frac 12 k^2 \left( \log k \right)^C \, .
\end{equation}
\end{lemma}
\medskip

Lemma~\ref{th:Tk} and \eqref{eq:tbuwlemTk} imply \eqref{eq:bound} and Proposition~\ref{th:bound} is proven.
\qed

\medskip

\begin{rem}
\label{rem:2}
The purpose of the assumption that  $\bbP(\logZ > \log 2)$ and $\bbP(\logZ < -\log 2 )$ are non zero is to guarantee that 
the  chain $X$ can exit $(-k,k)$. Under this assumption in fact $\bbP(X_{n+1}\ge k\vert X_n=x)>0$
for  
$x$ sufficiently close to $k$ (and analogous statement at $-k$): note in fact that $h_k(\pm k)= \pm (k-\log 2+ \log(1+e^{-2k})$.  If either $\bbP(\logZ > \log 2)=0$ or   $\bbP(\logZ < -\log 2 )=0$, the definitions of the chain $X^\approx$ should be modified by stipulating that the chain may step 
 to the absorbing state $\cC$ only  when $\vert x \vert \ge k-c$, with $c>0$ such that 
$(k-c)-h_k(k-c)(>0)$ and $(k-c)+h_k(-(k-c))(< 0)$ are in the interior of the support of $\gz$. Of course  in this case the probability of 
jumping to $\cC$ is no longer $1/2$: since 
\begin{equation}
\label{eq:1/2tcasece }
h'_k(\pm (k -c))\, =\,
1- \frac 1{e^c+1}- \frac 1{e^{2k-c}+1} \, <\,  1- \frac 1{e^c+1}\, ,
\end{equation}
we can choose this probability equal to $1- 1/(e^c+1)$. These changes do not affect Proposition~\ref{eq:Tk} because they only affect the constants $C$ and $k_0$, whose precise values are unimportant.
\end{rem}

\medskip


 \medskip
 
\medskip

\noindent
\emph{Proof of Lemma~\ref{th:Tk}.}
We recall that $\gz$ is centered and $\gs ^2=\int x^2 \gz(x) \dd x > 0$.
We claim that it suffices to show that for $k\ge k_0$
\begin{equation}
\label{eq:claimTk}
\inf_{x\in \bbR} \bbP_x\left( \tau_k \le k^2\right) \, \ge \, 2( \log k) ^{-C}=: p_k\, .
\end{equation}
In fact \eqref{eq:claimTk} implies that for $j\in \bbN$ and every $x$
\begin{equation}
\bbP_x\left( \tau_k >(j+1) k^2\right) \, \le \, \bbP_x\left( \tau_k > j k^2\right) \left(1- p_k\right)\, ,
\end{equation}
and therefore (also for $j=0$)
\begin{equation}
\bbP_x\left( \left\lceil\frac{\tau_k}{k^2} \right\rceil>j \right)\, =\, 
\bbP_x\left( \tau_k >j k^2\right) \, \le \,  \left(1- p_k\right)^j\, ,
\end{equation}
which implies 
\begin{multline}
\bbE_x \left[ \tau_k\right] \, \le \,  k^2 \bbE_x \left[ \left \lceil\frac{\tau_k}{k^2}\right \rceil\right] 
\, = \,  k^2 \sum_{j=0}^\infty \bbP_x\left( \tau_k >j k^2\right) \\ \le\,
k^2 \sum_{j=0}^\infty  \left(1- p_k\right)^j\, =\, \frac{k^2}{p_k} \, =\,  \frac 12 k^2 \left( \log k \right)^C\, .
\end{multline}

\medskip

Let us prove \eqref{eq:claimTk}:
it suffices to prove that both 
\begin{equation} 
\inf_{x\ge 0 } \bbP_x\left( \tau_k \le k^2, \, 
X_{\tau_k} >0\right)  \ \text{ and } \
 \inf_{x\le 0 } \bbP_x\left( \tau_k \le k^2, \, 
X_{\tau_k} < 0\right)
\end{equation} 
are bounded below by $2( \log k) ^{-C}$.  These two estimates are  \emph{equivalent} up to the fact 
that $\gz $ is not symmetric (but this affects only the choice of the constants and in an obvious way). So we just treat  the first expression, i.e.\ when the chain exits on the right.

We take this occasion also to remark that if we take two chains $X$ and $X'$, with $X_0=x$ and  $X'_0=x'$, both satisfying \eqref{eq:XMC}, coupled by the common randomness $( \logZ_n)$, we have
$X_{n+1}-X'_{n+1}= h_k(X_n)- h_k(X'_n)$, with $h_k(\cdot) $ increasing,  so the sign of $X_{n}-X'_{n}$ does not depend on $n$. 
Therefore 
$\bbP_0\left( \tau_k \le k^2, \, 
X_{\tau_k} >0\right)\ge \bbP_{x'}\left( \tau_k \le k^2, \, 
X_{\tau_k} >0\right)$  for $x \ge x'$. So it suffices  to check that 
\begin{equation} 
\label{eq:tbv}
\bbP_0\left( \tau_k \le k^2, \, 
X_{\tau_k} >0\right)\, \ge \,  2( \log k) ^{-C}\,.
\end{equation}
Moreover, what the process does outside $(-k, k)$ is irrelevant, so, when $\vert X_n\vert \ge k$, we replace $h_k(X_n)$ in 
\eqref{eq:XMC}  with $X_n$. This means that $h_k(x)$ is redefined as $h_k(x) \ind_{(-k,k)}(x)+ x\ind_{\bbR \setminus(-k,k)}(x)$.
Note that this corresponds to choosing a different increasing function $h_k(\cdot)$. 

To show that  \eqref{eq:tbv} holds we choose 
$\gd>0$ so that $\inf_{x\in[0,k)}\bbP_x( X_1-X_0 \ge \gd \vert \, X_0 =x) =p_\gd>0$. Note that this infimum 
is achieved at $x=k$ and, since we have trivialized the drift outside $(-k,k)$, the infimum over $x\in[0,k)$ can be taken over $x \in [0, \infty)$.   We introduce also the stopping time
\begin{equation}
\tau'_k\, :=\, \inf \left\{ n\ge \tau_{k- 2 \log k}: \, \vert X_n-(k- 2 \log k) \vert \ge 2 \log k - 2 \log \log k \right\}\, ,
\end{equation}
and we 
remark that
\begin{multline}
\label{eq:forclaimTk1}
\bbP_0\left( \tau_k \le k^2\, , X_{\tau_k} >0\right) \, \ge \, 
\bbP_0\Big( \tau_{k- 2 \log k} \le  k^2 /2 \textrm { and } X_{ \tau_{k- 2 \log k}} \ge  k- 2 \log k
\, , \\
\tau'_k-\tau_{k-2 \log k} \le (2 \log k)^2 \text{ and } X_{\tau'_k} \ge k-2 \log \log k
\, ,
\\
X_{\tau_{k- 2 \log \log k}+j}- X_{\tau_{k- 2 \log \log k +j-1} }\ge \gd \textrm{ for } j=1, 2, \ldots, \lceil (2/ \gd) \log \log k\rceil
\Big)\, .
\end {multline}
Let us read  this  complicated-looking expression: each item in the list that follows corresponds to the three lines 
of the right-hand side of \eqref{eq:forclaimTk1}.
\begin{enumerate}
\item  We ask that in less then $k^2/2$ steps the chain exits the interval $(-k+2 \log k,k-2 \log k)$ from the right: we will show that, by Brownian motion approximation,
this probability is positive uniformly in $k$. 
For this to work we need the drift to be negligible on the relevant (diffusive) time scale, which given the form of $h_k$ requires us to stop sufficiently far  (i.e.\ $2 \log k$) from the boundary.
\item 
We ask that  the chain exits on the right the interval of length $2(2 \log k - 2 \log \log k)$
centered in $k- 2 \log k$ in  $4 (\log k)^2$ steps or less. 
Note that in the unlikely (favorable) event that  $\tau'_k=\tau_{k-2 \log k}$ then 
$X_{\tau_{k- 2 \log  k}}\ge k -2 \log \log k$  is verified, so we can assume  
$X_{\tau_{k- 2 \log} }\in [k-2\log k,  k-2 \log \log k)$ and 
again we can treat this step by Brownian motion approximation. The drift is larger, but still negligible because the time span is much shorter than what we considered before. This will cost a probability factor bounded away from zero, like for (1).
We also remark that, on this event, $\tau'_k=\tau_{k-2 \log \log k}$.
\item  Once the chain is at a distance at most $2 \log \log k$ from the right boundary point $k$, we ask that it goes straight out of the domain by making steps of length at least $\gd>0$: 
this will cost a factor smaller than one (but bounded below by $p_\gd>0$) at each step. Since the number of steps is proportional to 
$\log \log k$ this costs a probability factor that vanishes like a power of $1/ \log k$. Note that we have modified the chain outside of $(-k,k)$, by removing the boundary repulsion, so that this ballistic strategy can be performed also once the chain is beyond $k$.   
\end{enumerate}

\medskip

Before starting the lower bound estimate for the right-hand side of
\eqref{eq:forclaimTk1}
we anticipate that, by the Strong Markov Property, 
we will be able to perform three separate estimates: each estimate  corresponding to one  of the three items in the list.
\smallskip 

Corresponding to the first item we aim at showing that 
\begin{equation}
\bbP_0\left( \tau_{k- 2 \log k} \le  k^2 /2  \textrm { and } X_{ \tau_{k- 2 \log k}} \ge  k- 2 \log k \right) 
\end{equation}
is bounded away from zero, uniformly in $k\ge k_0$. 
The event does not change if we replace $h_k(x)$ with $x$ for $\vert x \vert \ge k -2 \log k$ in defining the Markov chain $X$: so we will do so. Moreover  
we
define  the process 
$X_{t, k}:= X_{t k^2}/k$ if $ tk^2$ is an integer and otherwise we define $X_{t, k}$ by affine interpolation
so the trajectory is continuous.  We are going to show, via a  diffusion approximation,   that
the sequence of processes $( X_{\cdot, k})_{k \in \bbN}$, with $X_{\cdot, k}$ a random element 
of $C^0([0, T); \bbR)$ (with $T>0$ arbitrary), converges in law  to a Brownian motion with variance $\gs ^2=\int x^2 \gz (x) \dd x$. This is  formally stated in the following Lemma.

\medskip

\begin{lemma}
\label{th:SV}
For every continuous and bounded function 
$F$ from $C^0([0, T); \bbR)$ to $\bbR$ we have that 
\begin{equation}
\lim_{k \to \infty} 
\left \vert \bbE_0 \left[ F\left( X_{\cdot,k}\right)\right] - E_0 [ F(  B^\gs_\cdot)] \right\vert \, =\, 0\,,
\end{equation}
where, under $P_0$, $B^\gs _\cdot= \gs B_\cdot $ with $B_\cdot$ a standard  Brownian motion.
\end{lemma}
\medskip

\noindent
\emph{Proof of Lemma~\ref{th:SV}.}
We apply the diffusion approximation procedure detailed in 
 \cite[pp. 266--272]{StroockVaradhan}.  By  \cite[Assumptions (2.4)-(2.5)-(2.6), Theorem~11.2.3]{StroockVaradhan}  
it sufficed to check three conditions:
\smallskip

\begin{enumerate}
\item The (vanishing) drift condition: 
\begin{equation}
\label{eq:SV-cond1}
\lim_{k \to \infty } \sup_{y:\, \vert y \vert \le k -2 \log k } k \left \vert \int (x-y) \gz (x-h_k(y)) \dd x \right \vert \, =\, 0\,.
\end{equation} 
Note that we can restrict to $\vert x \vert \le k -2 \log k$ because the process has been modified so  to be centered outside of this interval.  \eqref{eq:SV-cond1} holds because, using first  $\int  x\gz  (x) \dd x=0$  and then
\eqref{eq:almost-id}, we obtain
\begin{equation}
\label{eq:SVintstep}
\left \vert \int (x-y) \gz (x-h_k(y)) \dd x\right \vert \, =\,\left \vert h_k(y)-y\right\vert \,\le\, \exp(-2 \log k)\, =\, \frac 1{k^2}\, , 
\end{equation} 
and \eqref{eq:SV-cond1} is verified.
\item The control of the variance:  
\begin{equation}
\lim_{k \to \infty } \sup_{y:\, \vert y \vert \le k -2 \log k }\left \vert \int (x-y)^2 \gz (x-h_k(y)) \dd x - \gs^2 \right \vert \, =\, 0\, , 
\end{equation}
and this is a direct consequence of the fact that, in the range of $y$ that we consider,  $\vert y-h_k(y)\vert $
is bounded by the boundary case $y= k -2 \log k$ and the resulting expression vanishes as $k \to \infty$.  
\item The control on large jumps: for every $\gep>0$
\begin{equation}
\lim_{k \to \infty } 
k \sup_{x\in \bbR} \bbP\left(\vert \logZ_1+ \left(h_k(x)-x\right) \ind_{[ -(k -2 \log k),+(k -2 \log k)]}(x) \vert > \gep k\right) \, =\, 0\, ,
\end{equation}
which is (largely) verified because $h_k(x)-x= O(1/k^2)$ for $\vert x \vert \le k-2 \log k$ and because $\logZ_1$ has finite exponential moments.
\end{enumerate}
\smallskip

This completes the proof of Lemma~\ref{th:SV}.
\qed

\medskip
 Lemma~\ref{th:SV} implies 
\begin{multline}
\bbP_0\left( \tau_{k- 2 \log k} \le  k^2 /2 \textrm { and } X_{ \tau_{k- 2 \log k}} \ge  k- 2 \log k\right) 
\stackrel{k \to \infty} {\longrightarrow}
 P_0 \left( \mathtt t _\gs \le 1/{2}, \,  B^\gs_{\mathtt t _\gs} =1\right)\, =: \, p_\gs \, ,
\end{multline} 
  and $\mathtt t _\gs$ is the hitting time
of $(-1,1)^\complement$ by $B^\gs _\cdot$. Hence for $k$ sufficiently large 
\begin{equation}
\bbP_0\left( \tau_{k- 2 \log k} \le  k^2 /2 \textrm { and } X_{ \tau_{k- 2 \log k}} \ge  k- 2 \log k\right) 
\, \ge \, \frac{p_\gs}2\, .
\end{equation}

\medskip

We now restart from time $ \tau_{k- 2 \log k}$ and use the Strong Markov Property. 
$X_{\tau_{k- 2 \log k}}\in [ k- 2 \log k, k-2 \log \log k)$  and we apply a Brownian motion 
approximation to the chain  $(X_{\tau_{k- 2 \log k+j}})_{j =0,1, \ldots}$ on the time scale $(\log k)^2$:
the steps are analogous to the ones of the previous step.
Also in this case it is 
  more concise to resort to the comparison argument explained right before \eqref{eq:tbv}, so that the starting point of our chain can and will be chosen equal to $k- 2 \log k$. 
  Therefore, by recentering,  i.e.\ by translating the system so that $k- 2 \log k$ becomes the origin (of course we have to shift accordingly the repulsion),  it suffices to show that 
\begin{equation}
\label{eq:thSV2}
\lim_{k \to \infty} 
\bbP_{0} \left( \tau_{ 2 \log k -2 \log \log k} \le (2 \log k)^2, \, X_{ \tau_{ 2 \log k -2 \log \log k} } >0 \right)\, =\,
P_0\left( \tau \le \mathtt t _\gs,
 B^\gs_{\mathtt t _\gs} =1\right)\,>0\,  ,
\end{equation}
But \eqref{eq:thSV2}
  is just  a close analog  of Lemma~\ref{th:SV}
 and the key point is that the proof is essentially the same up to replacing $k$ with $\log k$: note notably that now the repulsion is much stronger, $O(\exp(-2 \log \log k))= O( 1/ (\log k)^2)$ uniformly in the interval, but this yields a negligible drift because time is $O(\log k)$. 
 
 Therefore we arrive at: there exists $p'_\gs>0$ and $k_0$ such that 
 \begin{multline}
\label{eq:forclaimTk1.1}
\inf_{k \ge k_0}
\bbP_0\Big( \tau_{k- 2 \log k} \le  k^2 /2 \textrm { and } X_{ \tau_{k- 2 \log k}} \ge  k- 2 \log k
\, , \\
\tau'_k-\tau_{k-2 \log k} \le (2 \log k)^2 \text{ and } X_{\tau'_k} \ge k-2 \log \log k
\Big)\, \ge \, p'_\gs
\end {multline}

 Now we apply again the Strong Markov Property using the fact that $X_{\tau_k-2 \log \log k}\ge k -2 \log \log k$
 and the last estimate is just a product estimate that leads to (recall \eqref{eq:forclaimTk1})
 \begin{equation}
 \bbP_0\left( \tau_k \le k^2\, , X_{\tau_k} \ge k\right) \, \ge\, 
 p'_\gs \, p_\gd ^{\lceil (2/ \gd) \log \log k\rceil}\, \ge \, 
 \frac 1{(\log k)^C}\, ,
 \end{equation}
 with $C=(3/ \gd) \log(1/p_\gd)$ and $k$ sufficiently large. 
 This completes the proof of Lemma~\ref{th:Tk}.
 \qed

 \appendix

\section{Complementary results on the auxiliary chains}
\label{sec:various}

The   analysis of the basic properties of the $X$ chain can of course be found for example in the first two chapters of \cite{cf:BL}: by basic properties we mean existence and uniqueness of an invariant probability, that follow from  irreducibility and positive recurrence. We choose not to discuss these issues in detail  because we do give below more details about the $Y$ chain that is a bit more delicate to deal with -- in particular, the invariant measure is not normalizable -- and because for the $Y$ chain we need a few  results that depend on our restricted framework. 
For $Y$ we are going to exploit  \cite{cf:MT} and we will in particularly give a Lyapunov function to show recurrence: 
for $X$   we have positive recurrence because of the rather 
evident  \emph{confinement properties}  that can be made explicit using for example the Lyapunov function  $x \mapsto (\vert x\vert -k)_+^2)$ and
  \cite[Th.~11.3.4]{cf:MT}. 

\medskip

\noindent  
\emph{Proof of Proposition~\ref{th:Y}.} We distinguish here  among the case in which the support of $\gz$ is bounded away from $-\infty$ and when it is not. 

In the first case it is straightforward  to see that if $-c= \inf\{x:\, \gz(x)>0\}<0$ then 
the process eventually enters $(-\log(\exp(c) -1), \infty)$ and does not leave this set: in fact 
$-\log(\exp(c) -1)$ is the fixed point of $y\mapsto h(y)-c$.
Moreover $Y$ is irreducible, more precisely $\psi$-irreducible in the terminology for example of \cite{cf:MT},
with $\psi$ any probability with a  density that is positive on $(-\log(\exp(c) -1), \infty)$ and zero on the complement: this is a direct consequence of the fact that $\logZ$ has a density, of \eqref{eq:iterY} and of the fact that $-c$ is the left edge of the support of the transition probability. 
Recurrence of $Y$ follows for example by applying the criterion in \cite[Th.~3.1]{cf:Lamperti60},  or we  can apply 
\cite[Th.~8.0.2]{cf:MT} with a Lyapunov function equal to  $\log ((x-M)_+ +1)$, with $M=M(\gs)>0$ suitably chosen. 

If the support of $\gz$ is not bounded below, then the process is still $\psi$-irreducible and this time $\psi$ is any probability with positive density. Recurrence can be established in the same way: the repulsion from the left is very strong. An explicit Lyapunov function is    $x^2 \ind_{(-\infty, 0) (x)} +\ind_{[0, \infty)}(x) \log ((x-M)_+ +1) $ \cite[Th.~8.0.2]{cf:MT} . 

 So in both cases
$\nu$ exists, it is $\gs$-finite and it is unique
 \cite[Th.~10.4.9]{cf:MT}). Of course 
 $\nu$ is characterized by \eqref{eq:fornu}.
In particular for every bounded Borel set $B$
\begin{equation}
\label{eq:fromfornu}
\nu (B)\, =\, \int_B \left(\int_ \bbR \gz (z-h(y)) \nu (\dd y) \right) \dd z\, ,
\end{equation}
and we use this formula to show that $\nu(B)< \infty$ for every bounded Borel set.
To show this we use  the fact that $\gz \ge \gep \ind_{(a, b)}$ for suitably chosen positive constants $\gep, a$ and $b$. Assume that there exists a bounded Borel set $B_1$ with $\nu(B_1)=\infty$. Then there exists $x_0\in B_1$ such that $\nu ((x_0-1/n, x_0+1/n))=\infty$ for every $n$ and this directly yields $\int_ \bbR \gz (z-h(y)) \nu (\dd y) =\infty$ for every $z$ such that $z+h(x_0)\in (a,b)$.
But in this case, by \eqref{eq:fromfornu}, $\nu(B)=\infty$ for every $B \subset (a-h(x_0),b-h(x_0))$ of positive Lebesgue measure, which is impossible because $\nu$ is $\gs$-additive.

To establish $\nu (\bbR) = \infty$
we suppose that $\nu (\bbR) < \infty$,  and we assume that $\nu$ is a probability. In this case we remark that the process defined  recursively by
 $
S_{n+1}:=\,  \logZ_{n+1} + S_n
$ and by $S_0:= Y_0$
satisfies $S_n \le Y_n$ for every $n$. Hence $\bbP(S_n \ge \sqrt{n})\le \bbP(Y_n \ge \sqrt{n})$.
But the Central Limit Theorem implies that the left-hand side converges as $n \to \infty$ to a positive number,
while the right-hand side vanishes in the same limit because $(Y_n)$ is tight by the Ergodic Theorem applied to our (supposedly) positive recurrent process. Hence $\nu (\bbR) =\infty$. 

For the left tail property 
we start by remarking that using \eqref{eq:fornu} with  $g=\ind_{(a, b+1)}$,  $a<b$, and 
restricting the integral in the right-hand side to $y \in (-\infty, \log(e-1))$, i.e.\ 
 $h(y)\in(0,1)$, and to $x \in (a,b)$ 
\begin{equation}
\label{eq:fornu-lb}
\nu((a,b+1))\, = \, \nu \times \gz\left( \left\{(y,z):\, a < z+h(y)< b+1 \right\}\right)
\, \ge \, \nu \left( \left\{y:\, h(y)\in (0,1) \right\}\right)\gz((a,b))\, ,
\end{equation}
that is
$\nu((a, b+1)) \ge \nu((-\infty, \log (e-1)))\gz((a, b)) $. By choosing $a$ and $b$ such that $\gz((a, b))>0$ we see that
 $\nu((-\infty, \log (e-1)))< \infty$.  Therefore $\nu((-\infty, x])< \infty$ for every $x$ and the proof is complete.
 \qed
 
\medskip

\noindent
\emph{Proof of Lemma~\ref{th:pre-bound}.} 
By setting $g=\ind_{(-\infty, x]}$ in \eqref{eq:fornu} we see that for every $x \in \bbR$
\begin{equation}
\label{eq:Fnu-eq}
F_\nu(x)\, =\, \int F_{\gz} (x-h(y)) \nu (\dd y)\, .
\end{equation}
Since both $h(\cdot)$ and $F_{\gz}(\cdot)$ are non decreasing we have that
$y \mapsto F_{\gz}(x-h(y))$ is non increasing. So, by  
\eqref{eq:Fnu-eq}, we have that for every $x$ and every $z$ 
\begin{equation}
\label{eq:intStep-3.1}
F_\nu(x)\, \ge \, \int_{-\infty} ^z F_{\gz}(x-h(y)) \nu(\dd y) \, \ge \, F_{\gz}(x-h(z)) F_\nu(z)\, .
\end{equation}
Set $-\tilde c:= \inf\{y:\, F_\zeta(y)>0\} \in [-\infty, 0)$: one extracts form \eqref{eq:iterY-0} that $ \inf\{y:\, F_\nu(y)>0\}= -\log (\exp(\tilde c)-1)$.
Now consider the sequence $(x_j)$, $x_0$ the value we have chosen for $F_\nu(x_0)=1$,
and, for $j\in \bbN$, $x_{j}= x_{j-1}+ \rho$ with $\rho>0$ such that 
$x_0-h(x_1)=x_0-h(x_0+ \rho)> - \tilde c$, 
so that $q:= F_{\gz} (x_0-h(x_1)) >0$. Such a choice of $\rho>0$ is possible because $x_0>  -\log (\exp(\tilde c)-1)$ so
\begin{equation}
x_0-h(x_0)\, =\,  -\log(1+ \exp(-x_0))\, >\, - \tilde c\,.
\end{equation}
 Note moreover that for $j\in \bbN$ we have
$x_j -h(x_{j+1})= -\rho -\log(1+ \exp(-x_{j+1})) \ge -\rho -\log(1+ \exp(-x_{1}))= x_0-h(x_1)$.
So $F_{\gz} (x_j-h(x_{j+1})) \ge q$ for every $j$.
Therefore from  \eqref{eq:intStep-3.1} we infer that 
\begin{equation}
{F_\nu(x_{j+1})}\, \le \, \frac 1 {q^j} \, .
\end{equation}
Since $F_\nu$ is non decreasing this yields the claim.
\qed

\section{About the Fourier transform of $\zeta$}
With our Laplace transform notation 
$\widehat \zeta (-it)=\gp_\logZ(t)$ is the characteristic function of $\logZ$, that is the Fourier transform of the law of  $\logZ$.
In this Appendix we show how our hypotheses
\eqref{eq:LR} and \eqref{eq:Holder} on the distribution $\zeta$ imply that the characteristic function satisfies the bound
\begin{equation}
\label{eq:Fourier1/2}
\int_\bbR \frac{\vert \gp_\logZ (t)\vert}{(1+\vert t \vert)^{1/2}} \dd t \, < \infty\, ,
\end{equation} 
and Lemma~\ref{th:Leonardo} directly yields \eqref{eq:Fourier1/2}, thanks to our hypotheses
\eqref{eq:LR} and \eqref{eq:Holder}. 
We use $\hat f(t):= \widehat f (-it)= \int_\bbR f(x) \exp(it x) \dd x$.

\medskip 

\begin{lemma}
\label{th:Leonardo}
Assume that there exists $\theta>0$ such that $\vert f(x)-f(y) \vert \le \vert x-y \vert ^\theta$ for every $x\neq y \in \bbR$ and that there exists
$c>0$ such that $\int_\bbR \vert f (x) \vert \exp(c \vert x \vert) \dd x < \infty$. Then
\begin{enumerate}
\item  $f(x)= O(\exp(-b \vert x \vert))$ with  $b= c\theta/(1-\theta)$;
\item 
$\int_\bbR ({ \vert \hat f (t)\vert }/{(1+\vert t \vert)^{1/2}}) \dd t < \infty.
$
\end{enumerate} 
\end{lemma}
\medskip

\noindent
\emph{Proof.} 
It suffices to consider the case $x \to \infty$. We proceed by contradiction: assume that  there exists a sequence of positive numbers $(x_n)$ with $\lim_n x_n=+ \infty$ and $\vert f(x_n) \vert \ge 2 \exp(-b x_n)$ for every $n$. 
Then for every $x \in [x_n, x_n+ \exp(-bx_n/ \theta)]$ we have 
$\vert f(x) \vert \ge \exp(-b x_n)$. Therefore
\begin{equation}
\int_{x_n}^{ x_n+ \exp(-bx_n/ \theta)} \vert f(x) \vert \exp(c x)  \dd x
\, \ge \, \exp(-bx_n/ \theta) \exp(-b x_n) \exp(c x_n) \, =\,1\, , 
\end{equation} 
which is impossible because $\int_\bbR \vert f (x) \vert \exp(c \vert x \vert) \dd x < \infty$. 
Therefore (1) is established.

For (2) we start by remarking that (1) implies that the continuous function $f$ is in $\bbL^p$ for every $p\ge 1$,
in particular for $p=2$. 
Letting 
$(\go_2(h))^2\, := \, \int_\bbR \left( f(x+h) - f(x)\right)^2 \dd x $,
for $h \in [-1/2,1/2]$ and $K>2$ we have
\begin{equation}
\begin{split}
(\go_2(h))^2\, 
&\le\, \int_{-K}^K \left( f(x+h) - f(x)\right)^2 \dd x + C \int_{K-h}^\infty \exp(-2b x) \dd x
\\
&\le \, 2K \vert h \vert ^{2\theta} + C' \exp(-2bK)\, , 
\end{split}
\end{equation}
so if we choose $K=(\theta/b) \log (1/ \vert h \vert)$ -- if  $(\theta/b) \log (2)\le 2$ we  just choose a smaller value for  $b$ -- we have that there exists $C>0$ such that for $\vert h \vert \le 1/2$.
\begin{equation}
\go_2(h) \, \le \, C \vert h \vert ^\theta \sqrt{ \log(1/ \vert h \vert)}\, ;
\end{equation}
also evidently
$\sup_h \go_2(h) \le \sqrt{2} \Vert f \Vert_2$. These estimates imply that for every $\ga\in(0, \theta)$
\begin{equation}
\int_ \bbR \frac{ (\go_2(h))^2}{\vert h \vert ^{1+2\ga}} \dd h\, < \, \infty\,,  
\end{equation}
which implies \cite[Ch.~5: (46) and paragraph after (47)]{cf:Stein} that
\begin{equation}
\int_\bbR { \big \vert  \hat f (t)\big\vert ^2}
{\left(1+ \vert t \vert \right)^{2\ga}}\dd t \, < \, \infty\, . 
\end{equation}
Since
\begin{equation}
\left(
\int_\bbR \frac{ \big\vert \hat f (t)\big\vert }{(1+\vert t \vert)^{1/2}} \dd t \right)^2
\, \le \, \int_\bbR { \big \vert  \hat f (t)\big\vert ^2}
{\left(1+ \vert t \vert \right)^{2\ga}}\dd t
\int_\bbR \frac{ 1 }{(1+\vert t \vert)^{1+2\ga}} \dd t \,,
\end{equation}
the proof of part (2) is complete.
\qed

\section*{Acknowledgments}

\begin{sloppypar}
	We are very grateful to Leonardo Colzani for the proof of Lemma~\ref{th:Leonardo} and to Bernard Derrida for several enlightening discussions.
	G.G.\ is partially supported by ANR–19–CE40–0023 (PERISTOCH).
	The work of R.L.G.\ was carried out in the Department of Mathematics and Physics of Roma Tre University (Italy) and supported by the European Research Council (ERC) under the European Union's Horizon 2020 research and innovation programme (ERC CoG UniCoSM, grant agreement No.\ 724939 and ERC StG MaMBoQ, grant agreement No.\ 802901). 
\end{sloppypar}

\end{document}